\magnification 1200
\input amssym.def
\input amssym.tex
\parindent = 40 pt
\parskip = 12 pt
\font \heading = cmbx10 at 12 true pt
 at 22 true pt
\font \medheading =cmbx7 at 16 true pt
 at 7 true pt
\def \R{{\bf R}}
\def \Z{{\bf Z}}

\centerline{\medheading Stability of Oscillatory Integral Asymptotics}
\centerline{\medheading  in Two Dimensions}
\rm
\line{}
\line{}
\centerline{\heading Michael Greenblatt}
\line{}
\centerline{December 18, 2011}
\baselineskip = 12 pt
\font \heading = cmbx10 at 14 true pt
\line{}
\line{}
\noindent{\bf 1. Definitions and Background Information}.

\vfootnote{}{This research was supported in part by NSF grant  DMS-1001070} Let $S(x,y)$ be a smooth function defined on
 a neighborhood of the origin, and let $\phi(x,y)$ be a
smooth cutoff function supported in a sufficiently small neighborhood of the origin. The main object we are looking at is 
the oscillatory integral
$$J_{S,\phi}(\lambda) = \int_{\R^2}e^{i\lambda S(x,y)}\phi(x,y)\,dx\,dy \eqno (1.1)$$
Such integrals come up in a variety of settings in analysis and mathematical physics.
For a given $S(x,y)$ we are interested in the asymptotic behavior of $|J_{S,\phi}(\lambda)|$ as 
$\lambda \rightarrow \infty$. As in the author's earlier paper [G1], in this paper we are focusing on how the asymptotic behavior is affected by perturbations of $S(x,y)$. We will always assume $\lambda > 0$, as multiplying $\lambda$ by $-1$ just causes $J_{S,\phi}(\lambda)$ to be replaced by its
complex conjugate.

If the Taylor series coefficients of $S(x,y)$ at the origin are all zero, then
one can readily show directly that for no $\delta > 0$  is $|J_{S,\phi}(\lambda)|$ bounded by 
$C\lambda^{-\delta}$, so we do not conisder such situations. 
We also may always assume that $S(0,0) = 0$; in the integrand $(1.1)$ one may always factor out a
$e^{i\lambda S(0,0)}$, effectively replacing $S(x,y)$ by $S(x,y) - S(0,0)$ and $|J_{S,\phi}(\lambda)|$ will not be affected. 
It also does no harm to assume that $\nabla S(0,0) = 0$. For if $\nabla S(0,0) \neq 0$, then in $(1.1)$ one can integrate by parts arbitrarily many times 
to show that $J_{S,\phi}(\lambda)$ decays faster than any $C_n\lambda^{-n}$, and the same is true if $S$ is replaced by
a small perturbation of $S$.

\noindent Thus throughout this paper we will assume that the following conditions hold.
$$S(0,0) = 0 \,\,\,\,\,\,\,  \,\,\,\,\,\,\, \,\,\,  \,\,\,\,\,\,\,\nabla S(0,0) = 0 \eqno (1.2)$$
By resolution of singularities (see [AGV] ch. 7 for details), for a given real-analytic nonconstant $S(x,y)$ satisfying  $(1.2)$,
there is a positive $\delta$ and an integer $p = 0$ or $1$ such that for any $\phi(x,y)$ with sufficiently small support,
as $\lambda \rightarrow \infty$ one has an asymptotic development
$$J_{S,\phi}(\lambda) = A_{S,\phi}\,\lambda^{-\delta}(\ln\lambda)^p + o\big(\lambda^{-\delta}(\ln \lambda)^p\big) \eqno (1.3)$$
This $(\delta,p)$ is optimal in the sense that $A_{S,\phi} \neq 0$ as long as $\phi(x,y)$ is nonnegative (or nonpositive) with
$\phi(0,0) \neq 0$. In
[IM2] it is shown that $(1.3)$ still holds for general smooth $S(x,y)$ so long as $S(x,y)$ can be put in Case 1 or Case 2 superadapted coordinates. (See Definitons 1.9 and 1.10 below.) They also show that in these situations $A_{S,\phi}$ is always a multiple of $\phi(0,0)$. In the appendix, we will give explicit formulas for such $A_{S,\phi}$ which we will use in the proof of Theorem 2.3. These will be 
the same formulas that were shown to hold in the real-analytic case in [G2].

In the remaining smooth finite type situations, namely when $S(x,y)$ has Case 3 superadapted coordinates, one does not necessarily have asymptotics $(1.3)$. However, in [G2] it is shown there is some $\delta > 0$ independent of $\phi$ such that for some constant $A_{S,\phi}$ one has 
$$  |J_{S,\phi}(\lambda)| \leq A_{S,\phi}\lambda^{-\delta}\eqno (1.4)$$
This $\delta$ is such that for any $\delta' > \delta$,  so long as $\phi(x,y)$ is nonnegative with $\phi(0,0) \neq 0$ there is no constant $B_{S,\phi,\delta'}$ such that asymptotically one has 
$$ |J_{S,\phi}(\lambda)| \leq B_{S,\phi,\delta'} \lambda^{-\delta'} \eqno (1.5)$$
In view of the above, we make the following definitions as in [G1].

\noindent {\bf Definition 1.1.} The {\it oscillation index} of $S(x,y)$ at the origin is the $\delta$ for which $(1.3)$ holds when $S(x,y)$ has Case 1 or 2 superadapted coordinates, or for which $(1.4)-(1.5)$ holds when $S(x,y)$ has Case 3 superadapted
coordinates.

\noindent By [G2] and [IM2], $p = 0$ in $(1.3)$ when $S(x,y)$ has Case 1 superadapted coordinates or Case 3 superadapted coordinates in the real-analytic case, and $p = 1$ if it has Case 2 superadapted
coordinates. In view of this and $(1.4)-(1.5)$, we make the following definition.

\noindent {\bf Definition 1.2.} The {\it multiplicity} of the oscillation index is defined to be $1$ if $S(x,y)$ has Case 2 superadapted
coordinates, and is defined to be zero otherwise.

\noindent {\bf Definition 1.3.} $S(x,y)$ is said to be of {\it oscillatory type} $(-\delta, p)$ if $S(x,y)$ has oscillation index 
$\delta$ with multiplicity $p$.

\noindent We order the oscillatory types of phase functions 
lexicographically; we say that $(-\delta_1, p_1) < (-\delta_2,p_2)$ if $-\delta_1 < -\delta_2$ or if
$-\delta_1 = -\delta_2$ and $p_1 < p_2$. Thus a smaller type corresponds to faster oscillatory integral decay. 

\noindent {\heading Newton polygons and adapted coordinates.}

\noindent We now provide some relevant definitions that will be used throughout the paper.

\noindent {\bf Definition 1.4.} Assume that $S(x,y)$ is of finite type at the origin, and let $S(x,y) = \sum_{a,b} s_{ab}x^ay^b$ denote the Taylor expansion of $S(x,y)$ at the origin. For any $(a,b)$ for which $s_{ab} \neq 0$, let $Q_{ab}$ be the 
quadrant $\{(x,y) \in \R^2: 
x \geq a, y \geq b \}$. Then the {\it Newton polygon} $N(S)$ of $S(x,y)$ is defined as
the convex hull of the union of all $Q_{ab}$.  

A Newton polygon consists of finitely many (possibly zero) bounded edges of negative slope
as well as an unbounded vertical ray and an unbounded horizontal ray. 

\noindent {\bf Definition 1.5.} The {\it Newton distance} $d(S)$ of $S(x,y)$ is defined to be 
$\inf \{x: (x,x) \in N(S)\}$.

\noindent The line $x = y$ comes up so frequently in this subject it has its own name.

\noindent {\bf Definition 1.6.} The {\it bisectrix} is the line with  \~equation $x = y$.

The connection between Newton polygons and the asymptotics of oscillatory integrals and sublevel set measures is as follows.
It was shown by Varchenko [V] in the real-analytic case, and by Ikromov-M\"uller [IM1] in the smooth case, that for any $S(x,y)$
one can do a real-analytic (resp. smooth) coordinate change fixing the origin of the form $(x,y) \rightarrow (x - \phi(y),y)$ or 
$(x,y) \rightarrow (x,y - \phi(x))$, such that in the new coordinates the oscillation index $\delta$ of $S(x,y)$ satisfies 
$\delta = {1  \over d(S)}$. The multiplicity is equal to $1$ if the 
coordinate change can be made such that afterwards the bisectrix intersects $N(S)$ at a vertex. 
The multiplicity is zero otherwise. Furthermore, in any
coordinate system with the same origin, $d(S) \leq {1 \over \delta}$. Thus the coordinate systems where $d(S)$ are maximal are important to a
discussion of asymptotics, and correspondingly we have the following definition.

\noindent {\bf Definition 1.7.} A local coordinate system near $(0,0)$ is called {\it adapted} if $d(S) = {1 \over \delta}$ in this
coordinate system. 

\noindent We now come to the definition of superadapted coordinate systems, which are a refinement of the notion of adapted coordinate systems.

\noindent {\bf Definition 1.8.} For a given edge or vertex $e$ of $N(S)$, let $S_e(x,y) = \sum_{(a,b) \in e} s_{ab}x^ay^b$ denote
the sum of all terms of the Taylor expansion of $S(x,y)$ on $e$. 

\noindent {\bf Definition 1.9.} A local coordinate system near the origin is called {\it superadapted}
if for any compact edge $e$ containing the intersection $(d(S),d(S))$ of the bisectrix with $N(S)$, the functions $S_e(1,y)$ and $S_e(-1,y)$ both do not have zeroes of order $d(S)$ or greater other than possibly $y = 0$.

In [G2] it is shown that any smooth function can be put in superadapted coordinates. Also, as one would expect by symmetry, one can use the functions $S_e(x,1)$ and $S_e(x,-1)$ in place of $S_e(1,y)$ and $S_e(-1,y)$. We select the $y$-variable here for definiteness. 

\noindent As in [G2], we divide functions that are in superadapted coordinates into three cases:

\noindent {\bf Definition 1.10.} A function $S(x,y)$ in superadapted coordinates is said to be in Case 1 if the bisectrix intersects the Newton polygon $N(S)$ in the interior of a compact edge. It is said to be in Case 2 if the bisectrix intersects $N(S)$ at a vertex, and in Case 3 if the bisectrix intersects $N(S)$ in the interior of the horizontal or vertical ray. 

\noindent {\bf 2. Theorems and a little history}

A notable work on the effect of phase perturbations on the two-dimensional oscillation index and its multiplicity was given by Karpushkin
in [K1]-[K2]. Here he proved estimates that were uniform in a strong way, for the real-analytic case. Specifically, Karpushkin's theorem is as follows. Let $D_r$ 
denote the open disk in $\R^2$ of radius $r$ centered at the origin, and $E_r$ the open disk in 
${\bf C}^2$ of radius $r$ centered at the origin. For a function $f(x,y)$ real-analytic on $D_r$, let
$\tilde{f}(z_1,z_2)$ denote the unique holomorphic extension of $f(x,y)$ to $E_r$. Then Karpushkin's 
theorem for oscillatory integrals is as follows.

\noindent {\bf Theorem:} ([K1]-[K2]) Suppose $S(x,y)$ is real-analytic on $D_r$ satisfying $(1.2)$ with oscillatory type
$(-\delta,p)$. Then there is an $\eta > 0$, an $s < r$, 
and a positive constant $C_S$ depending on $S(x,y)$ such that if $f(x,y)$ is real-analytic on $D_s$
and $\tilde{f}(z_1,z_2)$ extends to a continuous function on $\bar{E}_s$ with $|\tilde{f}(z_1,z_2)| 
< \eta$ for all $(z_1,z_2) \in \bar{E}_s$ then for all $ \lambda > 0$ and all $\phi \in C_c^{\infty}(D_s)$ one has the estimate
$$|J_{S + f,\phi}(\lambda)| \leq C_{S,\phi}(1 +  \lambda) ^{-\delta} \ln(1 + \lambda)^p \eqno (2.1)$$
The proofs in [K1]-[K2] as well as in his other papers used
ideas from singularity theory, notably the theory of versal deformations, which is a way of converting arbitrary real analytic
perturbations of $S(x,y)$ into a finite list of canonical forms which then can be analyzed individually.

Another method for dealing with stability of  oscillatory integrals and sublevel set measures was introduced in [PSSt], where a 
slightly weaker version of Karpushkin's theorem is proved. [PSSt] uses a method which is sometimes referred to as the method 
of algebraic estimates. Some higher-dimensional theorems are also proved in [PSSt]; an example of Varchenko in [V]
shows the full analogues are not feasible. It should be mentioned that for the case of linear perturbations of smooth functions in
two dimensions, stability is proven in [IKeM]. In addition, there is an extensive body of research concerning analogous stability issues for complex analytic functions in several variables. 

 Our first theorem describes how the oscillatory type decreases or remains the same in a given direction, modulo finitely values, and gives criteria on a direction for the type to remain the same.

\noindent {\bf Theorem 2.1.} Suppose a smooth $S(x,y)$ satisfies $(1.2)$ and is in superadapted coordinates. Let $f(x,y)$ 
be a smooth function with nonvanishing Taylor expansion at $(0,0)$ and satisfyng $(1.2)$. Let $(-{\delta \over 2}, p)$ denote
 the oscillatory type of $S(x,y)^2 + f(x,y)^2$. There is a finite subset $I \subset
\R$ such that if $t \notin I$ then the oscillatory type of $S(x,y) + t f(x,y)$ is $(-\delta, p)$. This type is less than or equal
to the oscillatory type of $S(x,y)$, with equality holding in the following situations.

\noindent {\bf a)} $S(x,y)$ is in Case 1 superadapted coordinates. If $x + my = c$ denotes the equation of the edge of $N(S)$ intersecting the bisectrix, then equality holds iff $N(f) \subset \{(x,y) \in \R^2: x + my \geq c\}$. 

\noindent {\bf b)}  $S(x,y)$ is in Case 2 superadapted coordinates. Then equality holds iff there are $\infty \geq m_1 > m_2 \geq 0$ and $c_1, c_2 > 0$
 such that the lines
$x + m_1y = c_1$ and $x + m_2y= c_2$ intersect at the vertex $(d(S),d(S))$ of $N(S)$ and such that $N(f) \subset 
\{(x,y): x + m_1y \geq c_1\} \cap \{(x,y): x + m_2y \geq c_2\}$. When $m_1 = \infty$, the first line is taken to be $y = 
d(S)$.

\noindent {\bf c)} $S(x,y)$ is in Case 3 superadapted coordinates. Let $y = d(S)$ or $x = d(S)$ denote the equation of 
the infinite ray of $N(S)$ 
intersecting the bisectrix. Then equality holds iff $N(f) \subset \{(x,y): y \geq d(S)\}$ or $N(f) \subset \{(x,y): x \geq d(S)\}$
respectively.  

Some motivation for Theorem 2.1, courtesy [DK] where they considered the complex analytic case, is as follows. Let ${\delta \over 2}$ denote the oscillatory index of $S(x,y)^2 + f(x,y)^2$
as above, and let $D_r$ denote the open disk of radius $r$ centered at the origin.
Then by [G2] for example, ${\delta \over 2}$ is the supremum of the $\mu$ for which $\int_{D_r} (S(x,y)^2 + f(x,y)^2)^{-\mu} dx\,dy$ is finite for small $r$ (this is true for any function). However, by the Cauchy-Schwarz inequality, for
 $v_1^2 + v_2^2 = 1$ one has
$$|v_1S(x,y) + v_2f(x,y)|^{-2\mu} \geq (S(x,y)^2 + f(x,y)^2)^{-\mu}\eqno (2.2)$$
Integrating $(2.2) $with respect to $x$ and $y$ over $D_r$ shows that the oscillatory index of $v_1S(x,y) + v_2f(x,y)$ is at
most $\delta$. Scaling this, we see that the same is true for any $S(x,y) + tf(x,y)$.

\noindent On the other hand, for any $\mu < \delta$ one has
$$\int_{v_1^2 + v_2^2 = 1} \bigg(\int_{D_r} |v_1 S(x,y) + v_2f(x,y)|^{-2\mu}\,dx\,dy\bigg) = $$
$$\int_{D_r}\bigg( \int_{v_1^2 + v_2^2 = 1}|v_1 S(x,y) + v_2f(x,y)|^{-2\mu}\bigg)dx\,dy$$
Rotating the inner integral and using that $\delta \leq 1$ this becomes
$$\int_{D_r}\bigg( \int_{v_1^2 + v_2^2 = 1}v_1^{-{\mu}}(S(x,y)^2 + f(x,y)^2)^{-{\mu}}\bigg)dx\,dy$$
$$= C_{\mu} \int_{D_r} (S(x,y)^2 + f(x,y)^2)^{-{\mu}}dx\,dy \eqno(2.3)$$
The right hand side of $(2.3)$ is finite for any $\mu < \delta$, so the oscillation index of $v_1S(x,y) + v_2f(x,y)$ is greater than or equal to $\delta$ for 
almost all $(v_1,v_2)$ on the circle. Scaling this, one gets that the oscillation index of $S(x,y) + tf(x,y)$ is greater than or equal to $\delta$ for 
almost all $t$ as well. In view of the above, we conclude the oscillation index of $S(x,y) + tf(x,y)$ is equal to $\delta$ for almost all
$t$.

Thus Theorem 2.1 shows the oscillation index of $S(x,y) + tf(x,y)$ is equal to to 
twice that that of $S(x,y)^2 + f(x,y)^2$ for not just almost all $t$, but for all but finitely $t$, and furthermore there is 
agreement with regards to the multiplicity. In addition, sufficient and necessary 
conditions on $f(x,y)$ are provided for the oscillatory type of the generic $S(x,y) + tf(x,y)$ to be the same as that of $S(x,y)$.
In the course of our arguments we will see that the exceptional values of $t$ for a given $f(x,y)$ are determined
"effectively", meaning that they are explicitly determined through a finite sequence of finding roots $r$ of polynomials of one 
variable and coordinate changes of the form $(x,y) \rightarrow (x, y + rx^m)$.  

The functions  where equality holds in Theorem 2.1 are important for the other theorems of this paper. Thus we make the 
following definition.

\noindent {\bf Definition 2.2.} For a given $S(x,y)$, a smooth function $f(x,y)$ is said to be a {\it good direction} if equality holds in 
Theorem 2.1 for $f(x,y)$.

Our second theorem provides a smooth analogue of the uniform estimates of Karpushkin's theorem, but only in good directions.
For the Case 1 and Case 2 statements, we will make use of the norm $|f|_{r,N}$ on perturbation functions defined by 
$$|f|_{r,N} = \sum_{0 \leq \alpha, \beta \leq N} \sup_{(x,y) \in D_r}|\partial_x^{\alpha}\partial_y^{\beta} f(x,y)| 
\eqno (2.4)$$
For the Case 3 statement, which concerns perturbations in good directions in the real-analytic case, we use the following
variant of the norm $|f|_{r,N}$. Without loss of generality, assume the bisectrix intersects the horizontal ray. Then a good 
direction $f(x,y)$ can be written as $f(x,y) = y^dF(x,y)$ for some real-analytic $F(x,y)$. We define the norm $\|f\|_{r,N}$ by
$$\|f\|_{r,N} = \sum_{0 \leq \alpha, \beta \leq N} \sup_{(x,y) \in D_r}|\partial_x^{\alpha}\partial_y^{\beta} F(x,y)| 
\eqno (2.5)$$
\noindent {\bf Theorem 2.2.} Suppose a smooth $S(x,y)$  satisfies $(1.2)$ and is in superadapted coordinates. Let $(-\delta, p)$ denote the oscillatory type of $S(x,y)$.

\noindent Suppose $S(x,y)$ is in Case 1 or Case 2 superadapted coordinates. Then there are $r, N > 0$ 
 such that for some constants $C$ and $\eta$ the following holds. For all $\lambda > 2$, all $\phi \in C_c^{\infty}(D_r)$, and all good directions $f$ with  $|f|_{r,N} < \eta$ we have
$$|J_{S + f,\phi}(\lambda)| < C\|\phi\|_{C^1(D_r)}\lambda^{-\delta}(\ln \lambda)^p \eqno (2.6)$$
If $S(x,y)$ is real-analytic and in Case 3 superadapted coordinates, then there are $r, N > 0$ 
such that for some constants $C$ and $\eta$ the following holds. For all $\lambda > 2$, all $\phi \in C_c^{\infty}(D_r)$, and all $f$ real-analytic on a neighborhood of $\bar{D}_r$ such that $f$ is a good direction with $\|f\|_{r,N} < \eta$ we have
$$|J_{S + f,\phi}(\lambda)|< C\|\phi\|_{C^1(D_r)} \lambda^{-\delta} \eqno (2.7)$$
{\bf Remark 1.} By Karpushkin's results, in the real-analytic situation $(2.6)-(2.7)$ hold in all directions if one uses the magnitude of
the complex extension of a perturbation function as its norm.

\noindent {\bf Remark 2.} For general smooth $S(x,y)$ in Case 3 superadapted coordinates, it can be shown the Case 1 and 2
statement $(2.6)$ still holds for $\phi(x,y)$ supported in some $D_s \subset D_r$, where $s$ depends on $f$. The author does not 
know if the full statement holds in the general smooth Case 3 situation.

Our third and final theorem shows that in Case 1 for all smooth phases, and in Case 3 for real-analytic functions, for some $\alpha > 0$ the coefficient $A_{S,\phi}$ in $(1.3)$ is locally 
an $\alpha$-Lipschitz function of $f$ in the norms $|f|_{r,N}$ and $\|f\|_{r,N}$ respectively, again when $f$ is restricted to good directions. (In other directions the coefficient abruptly drops to zero.)

\noindent  {\bf Theorem 2.3.}   If $S(x,y)$ is in Case 1 of superadapted coordinates and $r$, $\eta$ and $N$ are as in Theorem 2.2, then there are $\alpha,  B > 0$ 
such that if $f_1$ and $f_2$ are good directions with  $|f_1|_{r,N}, |f_2|_{r,N} < \eta$, then the leading coefficient of the asymptotics $(1.3)$ satisfies
$$|A_{S + f_2,\phi} - A_{S + f_1\phi}| \leq B|f_2 - f_1|_{r,N}^{\alpha}|\phi(0,0)| \eqno (2.8)$$
\noindent Suppose $S(x,y)$ is in Case 3 of superadapted coordinates and $r$, $\eta$ and $N$ are as in Theorem 2.2. There 
are $\alpha,B > 0$ such that if 
$f_1(x,y)$ and $f_2(x,y)$ are good directions real-analytic on a
neighborhood of $\bar{D}_r$ with  $\|f_1\|_{r,N}, \|f_2\|_{r,N} < \eta$, then
$$|A_{S + f_2,\phi} - A_{S + f_1,\phi}| \leq B\|f_2 - f_1\|_{r,N}^{\alpha}\|\phi\|_{L^{\infty}} \eqno (2.9)$$
{\bf Remark 1.} The Case 2 analogue of $(2.8)-(2.9)$ is actually false as simple examples  show; take $S(x,y) = x^dy^d$ and then perturb by $f_1(x,y) = \epsilon x^dy^d$ and $f_2(x,y) = \epsilon x^ey^f$ where $e < d$ and $f > d$.
However, if one restricts to $f(x,y)$ such that $N(f) \subset N(S)$, then the analogue does hold; this is an immediate consequence
of the formulas in [G2] and the appendix of this paper. 

\noindent {\bf Remark 2.} There is no smooth Case 3 version of $(2.9)$ since one does not 
necessarily have an asymptotic expansion for a smooth Case 3 phase. 

\noindent {\bf Remark 3.} The proof of Theorem 2.3 will reveal that one can take $f_1$  or $f_2$ to be the zero function in $(2.8)$ and $(2.9)$.

\noindent {\bf Remark 4.} One also has an analogue of Theorem 2.3 where $\phi$ is also allowed to change; one just bounds $|A_{S + f, \phi_2 - \phi_1}|$ using Theorem 2.2 and adds the result to $(2.8)$ or $(2.9)$.

\noindent {\bf 3. Preliminary lemmas and the proof of Theorem 2.1.}

\noindent  For a polynomial $f(x)$, let $ord_x(f)$ denote the order of a zero of $f$ at $x$, with $ord_x(f) = 0$ if $f(x) \neq 0$. We start by proving a lemma concerning polynomials in one variable.

\noindent {\bf Lemma 3.1.} Suppose $p(x)$ and $q(x)$ are polynomials in one variable, neither identically zero. Let  let $\R_0$ denote $\R - \{0\}$, and let $m = \sup_{x \in \R_0}\min(ord_x(p),ord_x(q))$. Then there is a finite $I \subset \R$ such that if $t \notin
I$, then $ord_x(p + tq) \leq \max(1, m)$ for all $x \in \R_0$.

\noindent {\bf Proof.} The result is immediate if $p(x)$ and $q(x)$ are multiples of each other, so we assume this is not the case.
 Let $J$ be the set $\{x \in \R:p(x)q'(x) - p'(x)q(x) = 0,$ $q(x) = 0$, or $p(x) = 0\}$. $J$ is finite since
$p(x)q'(x) - p'(x)q(x)$ can't be identically zero; if it were by the quotient rule for derivatives $({p \over q})'(x)$ would be zero
 for all $x$, contradicting the assumption that $p(x)$ and $q(x)$ are not multiples of each other.

 I claim that for all $t \in \R$ the
 function $p(x) + tq(x)$ has no zero of order greater than 1 on $\R - J$. For suppose $p(x) + tq(x)$ had a zero of order at least 2 at some $x_0 \in \R - J $. Then $p(x_0)q'(x_0) = -tq(x_0)q'(x_0)= p'(x_0)q(x_0)$. So $p(x_0)q'(x_0) - p'(x_0)q(x_0) = 0$, contradicting that $x_0 \in \R - J$.

Now suppose $x_0$ is in the finite set $J$. Since Lemma 3.1 and our proof so far are symmetric in $p$ and $q$, without loss of generality we may assume
that $ord_{x_0}(p) \leq ord_{x_0}(q)$ (which may be zero.) Write $p(x) = (x - x_0)^{ord_{x_0}(p)}r(x)$ and $q(x) = (x - x_0)^{ord_{x_0}(q)}s(x)$. Thus we have
$$p(x) + tq(x) = (x - x_0)^{ord_{x_0}(p)}(r(x) + t (x - x_0)^{ord_{x_0}(q) - ord_{x_0}(p)}s(x))\eqno (3.1)$$
Let $K$ be a closed interval centered at $x_0$ such that $r(x) \neq 0$ on $K$. Then letting $u(x) = {s(x) \over r(x)}$, on $K$ we
may write
$$p(x) + tq(x) = r(x)(x - x_0)^{ord_{x_0}(p)}(1 + t (x - x_0)^{ord_{x_0}(q) - ord_{x_0}(p)}u(x))\eqno (3.2)$$
Here $u(x_0) \neq 0$. Consider the case that $ord_{x_0}(q) > ord_{x_0}(p)$.  Note that  there is an
interval $K$ centered at $x_0$ such that for any $t$ the factor $1 + t (x - x_0)^{ord_{x_0}(q) - ord_{x_0}(p)}u(x)$ has at most
two zeroes on $K$, each with multiplicity one, and $x = x_0$ cannot be a zero. Thus $p(x) + tq(x)$ has a zero of 
order $ord_{x_0}(p)$ at $x = x_0$, and no other possible zeroes of multiplicity greater than 1 on $K$. This is better than
what we need.

Now consider the case where  $ord_{x_0}(q) = ord_{x_0}(p)$. Then we write $u(x) = u_0 + (x - x_0)^mv(x)$ for some
 $m > 0$ and some $v(x)$ with $v(x_0) \neq 0$ (again we use that $p$ and $q$ are not multiples of each other). Now we have
$$p(x) + tq(x) = r(x)(x - x_0)^{ord_{x_0}(p)}(1 + t u_0 + t(x - x_0)^mv(x))$$
If we exclude $t = -{1 \over u_0}$, then exactly as in the $ord_{x_0}(q) > ord_{x_0}(p)$ situation there is some $K$ 
centered at $x_0$ such that $p(x) + tq(x)$ has a zero of order $ord_{x_0}(p)$ at $x = x_0$, and no other possible zeroes 
of multiplicity greater than 1 on $K$. Thus regardless of whether or not  $ord_{x_0}(q) > ord_{x_0}(p)$, excluding one
possible value of $t$ if necessary we can find an interval $K$ centered at $x_0$ 
such that the zeroes of $p(x) + tq(x)$ on $K$ are of order at
most $\max(1, ord_{x_0}(p))$, which is at most $\max(1,m)$ as long as $x_0 \neq 0$. Since $J$ is finite, we can cover
$\R_0 \cap J$ with finitely many intervals $K_i$ with this property. Since all zeroes of all $p(x) + tq(x)$ are of order 1 on 
$\R_0 - J$ we are done with the proof of the lemma.

Recall that for a smooth function $g(x,y)$ satisfying $(1.2)$, being in superadapted coordinates means that
if $e$ is a compact edge of $N(g)$ intersecting the bisectrix then the polynomials $g_e(1,y)$ and $g_e(-1,y)$ do not have zeroes of order $d(g)$ or higher, except possibly at $y = 0$. In this situation, by [G2] we necessarily have $|J_{g,\phi}| \leq C_{g,\phi}\lambda^{-\delta}
\ln(\lambda)^p$ where $(-\delta,p)$ is the type of $g$ at the origin (even in Case 3 superadapted coordinates where asymptotics
may not exist). 

Suppose we consider functions $g(x,y)$ defined on
 the semidisk $D_r \cap \{x > 0\}$ that are of the form $s(x^{1 \over n},y)$ for some
smooth $s(x,y)$ on a neighborhood of the origin. Then we may form the Newton polygon of $g(x,y)$ as before, except instead of having exponents
in $\Z \times \Z$, we consider $g(x,y)$'s Taylor expansion as having exponents in ${1 \over n}\Z \times \Z$. In this paper we 
will be only considering $g(x,y)$ whose Taylor expansions only contain terms of total degree two or more, so we will
 assume that throughout. Note that this implies that $d(g) \geq 1$.

 For $\lambda > 0$ and
 a smooth compactly supported function $\phi(x,y)$ on  $D_r \cap \{x \geq 0\}$ we look at the oscillatory integral
$$J_{g,\phi}^+(\lambda) = \int_{\{x > 0\}}e^{i\lambda g(x,y)}\phi(x,y)\,dx\,dy\eqno (3.3)$$
We define the type of $g$ to be the minimal $(-\delta,p)$ such that  one has an estimate  $|J_{g,\phi}^+(\lambda) | \leq C_{g,\phi}\lambda^{-{\delta}}\ln(\lambda)^p$ valid for all $\phi(x,y)$ with sufficiently small support. 

\noindent {\bf Lemma 3.2.}
Suppose the leftmost vertex of $N(g)$ is of the form $x^ay^b$ for integers $a$ and $b$ 
and that one of the following two situations holds.

\noindent 1) $N(g)$ does not intersect the bisectrix at a vertex and for each compact edge
$e$ of $N(g)$ any zeroes of $g_e(1,y)$ are of order $< d(g)$, or of order 1 in the case that $d(g) = 1$.

\noindent 2) $N(g)$ intersects the bisectrix at a vertex and for each compact edge $e$ of $N(g)$ the 
zeroes of $g_e(1,y)$ are of order $\leq d(g)$.

\noindent Then  $g$ is of type $(-{1 \over d(g)},p)$, where $p = 1$ if $N(g)$ intersects the bisectrix at a vertex and $p=0$
otherwise.  Furthermore for sufficiently small $r > 0$ and $0 < \epsilon < {1 \over 2}$ there are constants $c_{g,r}$ and $c_{g,r}'$ such that 
$$c_{g,r}\epsilon^{-{1 \over d(g)}}\ln(\epsilon)^p \leq |\{(x,y) \in D_r^+: |g(x,y)| < \epsilon\}| \leq c_{g,r}'\epsilon^{-{1 \over d(g)}}\ln(\epsilon)^{p'} \eqno (3.4)$$
Here $D_r^+$ are the points in $D_r$ for which $x > 0$, and $p = p'$ unless $d(N) = 1$ and some $g_e(1,y)$ has a zero, in which case $p= 0$ and $p' = 1$.

\noindent {\bf Proof.} We will not prove this lemma in detail since very similar arguments appear in [G2]; related arguments are
also in section 4 of this paper. Briefly, the lemma is true for the following reasons.
One may split up a neighborhood of the origin into "slivers" touching the origin and each sliver is a
 subset of $\{x > 0\}$ or $\{x < 0\}$. One finds upper bounds for 
$J_{g,\phi}^+(\lambda) $ or $|\{(x,y) \in D_r^+: |g(x,y)| < \epsilon\}|$  by adding up estimates corresponding the $\{x > 0\}$ slivers, and the estimate for a given sliver is given
by doing Van der Corput-type arguments in $x$ or $y$ and then integrating the result in the other variable.
These Van der Corput arguments
can almost always be made in the $y$ direction or for a first derivative in the $x$-direction, and these arguments carry over 
directly to the current setting. The sole exception is for the sliver corresponding to the left-most
vertex, so we require that vertex to be of the form $x^ay^b$ for integers $a$ and $b$. Then the arguments of [G2] carry 
over for that sliver as well. 

Sharpness of the exponents obtained in this way comes from directly estimating $J_{g,\phi}^+(\lambda)$ from below on certain rectangles determined by $N(g)$, for $\phi$ such that $\phi(0,0) \neq 0$. Similarly, the lower bounds of $(3.4)$ come from 
evaluating the measure of such rectangles arising from $N(g)$ on which $|g(x,y)| < \epsilon$. This concludes our
discussion of the proof.

Now suppose $g(x,y)$ and $h(x,y)$ are two functions of the form $s(x^{1 \over n},y)$, $s$ smooth on a neighborhood of 
the origin, such that the exponents of their left-most vertices are integers. We direct our attention to $N(g + th)$ for 
various $t \in \R$. There is a finite $I \subset \R$  such that if $t \notin I$, then the vertices of $N(g + th)$ are all vertices of
$N(g)$ or $N(h)$; one excludes the $t$ for which a term in $g(x,y)$'s Taylor expansion corresponding a vertex
of $N(g)$ is exactly canceled out by a term in $th(x,y)$'s Taylor expansion, or vice versa. In fact there is a single
Newton polygon $N$ which is equal to $N(g + th)$ for all $t \notin I$. 

Note that if $t \notin I$ and $g + th$ satisfies the hypotheses of 
of Lemma 3.2, then $g + th$ is of type $(-{1 \over d(N)},p)$ and equation $(3.4)$ holds for $g + th$.
For the purposes of proving Theorem 2.1,
we are interested in knowing for which $t$ does this happen. For this we have the following lemma.

\noindent {\bf Lemma 3.3.} For a given compact edge $e$ of $N$ let $g^e(x,y)$ and $h^e(x,y)$ denote the sum of the terms of the Taylor expansions of $g$ and $h$ respectively whose exponents are on $e$. As above, let $N$ be the common Newton polygon
of $g + th$ for all but finitely many $t$.

\noindent {\bf a)} Suppose $N$ does not intersect the bisectrix at a vertex and for each compact edge $e$ of $N$
$ \sup_{y_0 \in \R_0} \min(ord_{y_0}(g^e(1,y)), ord_{y_0}(h^e(1,y))) < d(N)$. Then there is a finite $J \subset \R$ such that if $t \notin J$, then any zero of any $(g + th)_e(1,y)$ besides $y = 0$ is of order less than $d(N)$ if $d(N) > 1$, and is of order $1$ if $d(N) = 1$.

\noindent {\bf b)} Suppose $N$ intersects the bisectrix at a vertex  and for each compact edge $e$ of $N$
$ \sup_{y_0 \in \R_0} \min(ord_{y_0}(g^e(1,y)), ord_{y_0}(h^e(1,y))) \leq d(N)$. Then there is a finite $J \subset \R$ such that if $t \notin J$, then any zero of any $(g + th)_e(1,y)$ other than $y = 0$ is of order less than or equal to $d(N)$.

\noindent {\bf c)} Let $(-{\delta \over 2},p)$ denote the type of $g(x,y)^2 + h(x,y)^2$.  In both part a) and part b),  if $t \notin J$ then $g(x,y) + th(x,y)$ is of type $(-\delta, p)$. Furthermore, $g(x,y)^2 + h(x,y)^2$ satisfies the hypotheses of Lemma 3.2, as does
$g(x,y) + th(x,y)$ if $t \notin J$.

\noindent {\bf Proof.} By Lemma 3.1,  in the setting of part a) if $d(N) > 1$ there is necessarily a finite $J \subset \R$ such that if $t \notin J$, then for each compact $e$ all zeroes of the polynomial $(g + th)_e(1,y)$ of order less than $ d(N)$, while if $d(N) = 1$ any zeroes
are of order 1. In the setting of part b) one can similarly remove a finite $J$ such that the zeroes are of order at most $d(N)$.
Thus by Lemma 3.2, in the settings of both part a) and part b) the type
of $g + th$ is $(-{1 \over d(N)}, p)$ where $p = 1$ iff $N$ intersects the bisectrix at a vertex. Furthermore, equation $(3.4)$ holds.

We now examine the function $g(x,y)^2 + h(x,y)^2$. Note that the vertices of
$N(g^2)$ are exactly $\{2v: v$ is a vertex of $N(g)\}$, and the edges of $N(g^2)$ are $\{2e: e$ is an edge of $N(g)\}$. (Here $2e$ denotes the dilate of $e$ by a factor of 2 and similarly for $2v$). Thus $N(g^2)$ is $2N(g)$, the dilate by 2 of $N(g)$. Furthermore, if $e$ is a compact edge of $N(g)$, then
$(g^2)_{2e}(1,y) = (g_e(1,y))^2$. The analogous statements also hold for $N(h^2)$. 

Thus the Newton polygon of $g^2 + h^2$
is exactly $2N$, and for any edge $e$ of $2N$, $(g^2 + h^2)_e(1,y) = (g^e(1,y))^2 + (h^e(1,y))^2$. Since $(g^e(1,y))^2$ 
and $(h^e(1,y))^2$ are nonnegative functions, the maximum order of any zero of  $(g^2 + h^2)_e(1,y)$ for $y \neq 0$ is
given by $\sup_{y_0 \neq 0} \min(ord_{y_0}((g^e(1,y))^2,(h^e(1,y))^2) = 2\sup_{y_0 \neq 0} \min(ord_{y_0}(g^e(1,y), h^e(1,y)))$. By assumption, this quantity is less than or equal to $2d(N) = d(2N)$ when $2N$ intersects the bisectrix at a
vertex, and is less than $d(2N)$ when $2N$ does not. Thus by Lemma 3.2, the type of $g^2 + h^2$ is given
by $(-{1 \over 2d(N)}, p)$ and $g^2 + h^2$ satisfies $(3.4)$ with $p = p'$ since $d(2N) > 1$. 
 Since $g + th$ is of type $(-{1 \over d(N)},p)$ for $t \notin J$, this completes the proof.

\noindent {\bf Proof of Theorem 2.1.} 

The argument will proceed as follows. We will consider $(x,y) \in D_r \cap \{ x> 0\}$, as analogous
arguments always work for $D_r \cap \{x < 0\}$. For a given $f$, we will either be able to apply Lemma 3.3 immediately, or
 through an iterative process we will construct a smooth $p(x)$ such that for some $n$,
if $S'(x,y) = S(x,y -p(x^{1 \over n}))$ and $f'(x,y) = S(x,y -p(x^{1 \over n}))$, then $S'(x,y)$ and $f'(x,y)$ will satisfy the 
conditions of Lemma 3.3. Thus in these situations, except for finitely many $t$ we will have $|J_{\phi,S' + tf'}(\lambda)^+| \leq C \lambda ^{-\delta}(\ln(\lambda))^p$, where $(-{\delta \over 2},p)$ denotes the type of $(S')^2 + (f')^2$, which by $(3.4)$ is the same as the type of $S^2 + f^2$  when viewed as a smooth function on all of $D_r$ (these coordinate changes do not affect the measures). Also, the coordinate changes are such that $J_{\phi,S' + tf'}^+(\lambda) = J_{\phi,S + tf}^+(\lambda)$, and if one adds 
the $x < 0$ and $x > 0$ portions of this, one gets $|J_{\phi,S + tf}(\lambda)| \leq C \lambda ^{-\delta}(\ln(\lambda))^p$.
That this type $(-\delta,p)$ is optimal follows from adding  over the $x < 0$ and $x > 0$  portions the lower bounds of $(3.4)$
and using that the oscillatory type of a smooth function $a(x,y)$ is equal to the supremum of the $(-\theta,q)$ for which one has an estimate $|\{(x,y) \in D_r: |a(x,y)| < \epsilon\}| \geq c_{g,r}'\epsilon^{\theta}\ln(\epsilon)^{q}$. (This was proved in [G2]; for the real-analytic case, see [AGV] Ch. 7.) 

In summary, to prove Theorem 2.1 it will suffice to show that either $S$ and $f$ satisfy the conditons of Lemma 3.3 to start with, or the derived $S'$ and $f'$ satisfy the conditions of Lemma 3.3. We will then identify when $S$ has the same type as the generic $S + tf$; this will determine whether or not $f$ is a good direction. It will transpire that good directions occur in some cases
when $S$ and $f$ initially satisfy the conditions of Lemma 3.3, and never when one performs the coordinate changes. The former
situations will correspond to the directions stipulated in Theorem 2.1.

We now begin the main argument. Let $N$ be the
polygon such that $N(S^2 + f^2) = 2N$. Let $I \subset \R_0$ such that if $t \in \R_0 - I$ then no vertices of $N(S)$ or $N(tf)$ get 
cancelled out in adding $S + tf$. So in particular, if $t\in \R_0 - I$, then as in the proof of Lemma 3.3 we have $2N(S + tf) = N(S^2 + f^2)$. If the bisectrix intersects $N$ at a vertex $v$ then either $N(f)$ or $N(S)$ has a vertex at $v$. In the first case each $f_e(1,y)$ has no zero of order greater
than $d(N)$ besides $y = 0$, and in the latter case the same holds for each $S_e(1,y)$.  Thus the conditions of Lemma 3.3 are
satisfied, with $f$ a good direction in the first case. If the bisectrix intersects $N$ in the interior of one of the infinite rays, then the bisectrix intersects either $N(f)$ or
$N(S)$ in the interior of this ray, and with the same Newton distance $d(N)$. In the former case each $f_e(1,y)$ has no zero of order $d(N)$ or higher, other than possibly $y = 0$, and in the latter case the same is true for each $S_e(1,y)$. Thus once again the hypotheses of Lemma 3.3 are immediately satisfied, with $f$ a good direction in the former case.
So from now on we may assume the bisectrix intersects $N$ in the interior of a compact edge which we denote by $e_0$. 

As before, let $S^{e_0}(x,y)$ and $f^{e_0}(x,y)$ respectively denote the sum of the terms of $S$'s and $f$'s Taylor expansion
that are on $e_0$; in either case there may be one, or even no terms. Suppose $S^{e_0}(x,y)$ has at least one term and each exponent $(a,b)$ appearing in 
$S^{e_0}(x,y)$ satisfies $b < a$. Then if $b_0$ denotes the maximum of such $b$, $S^{e_0}(1,y)$ is a polynomial of 
degree $b_0 < d(N)$. Hence the zeroes of  $S^{e_0}(1,y)$ are of order at most $b_0$, and thus  $\min_{y_0 \in \R_0}(ord_{y_0}(S^{e_0}(1,y)), ord_{y_0}(f^{e_0}(1,y))) < d(N)$. Furthermore, if $e'$ is an edge of $N$ below the bisectrix, a nonzero $f^{e'}(1,y)$, $S^{e'}(1,y)$ has 
degree $< d(N)$, while if $e'$ is above the bisectrix then a nonzero $f^{e'}(x,1)$ or $S^{e'}(x,1)$ has degree $< d(N)$. Thus these polynomials have no zeroes of degree $d(N)$ or larger. By mixed homogeneity of $f^{e'}(x,y)$
and $S^{e'}(x,y)$, in the latter two
cases the same must be true for $f^{e'}(1,y)$ or $S^{e'}(1,y)$, other than $y = 0$. Thus the conditions of Lemma 3.3 are satisfied. Furthermore, since $S$ is in superadapted coordinates and the bisectrix intersects $N(S)$ above the line containing $e_0$, the type of $S$ is greater than the type of the generic $S + tf$ and hence $f$ is not a good direction.

The argument of the previous paragraph still applies if each exponent $(a,b)$ satisfies $b > a$ instead of 
$a < b$, since one can reverse the roles of the $x$ and $y$ variables in this situation. We can also replace $S^{e_0}(x,y)$ by
$f^{e_0}(x,y)$  and the above argument still gives that Lemma 3.3 applies. This time the generic $S + tf$ is of the same type as $S$  iff $S^{e_0}(x,y)$ has terms both above and below the bisectrix. So $f$ is a good direction when this happens.

Thus we may restrict our attention from now on to the situation where the bisectrix intersects $N$ in the interior of a compact 
edge $e_0$, and where
$f^{e_0}(x,y)$ are $S^{e_0}(x,y)$ are either the zero polynomial, or are nonzero and have exponents that aren't all strictly 
above or strictly below the bisectrix. Suppose first $S^{e_0}(x,y)$ is nonzero. Since its exponents aren't all strictly above or 
below the bisectrix, either $N(S) \cap e_0 = \{(d(N),d(N))\}$, or $N(S)$ has a compact edge $e$
intersecting the bisectrix such that $e \subset e_0$. In the former case, $S^{e_0}(1,y)$ has no zeroes at all other
than $y = 0$, and in the latter case since $S$ is in superadapted coordinates $S^{e}(1,y)$ has no zeroes of degree $d(N)$ or
more other than possibly $y = 0$. Any $f^{e'}(1,y)$ or $S^{e'}(1,y)$ for an edge $e'$ of $N$ not intersecting the bisectrix will have zeroes of order less than $d(N)$ similarly to two paragraphs ago. Thus the conditions of Lemma 3.3
are satisfied. If $N(S)$ has a vertex on the bisectrix, then $S$ of type $(-{1 \over d(N)}, 1)$ while the generic $S + tf$ is of type
$(-{1 \over d(N)},0)$. Hence $f$ is not a good direction. Otherwise, $S$ is of type $(-{1 \over d(N)},0)$ so, $f$ is a good direction.

Thus it remains to consider the case where $S^{e_0}(x,y)$ is the zero polynomial, so that  $(S + tf)_{e_0}(x,y) = tf_{e_0}(x,y)$. 
If each zero of $f_{e_0}(1,y)$ is of
order $< d(f) = d(N)$, then $S + tf$ is in superadapted coordinates and like in the above cases Lemma 3.3 applies. Note that 
$f$ is not a good direction here as $d(S) > d(S + tf)$. So we may
devote our attention to where $S^{e_0}(x,y)$ is the zero polynomial and $f_{e_0}(1,y)$ has a zero of order at least $d(f) = d(N)$; in essence this situation is the crux of the proof.

Let the
equation of the line containing $e_0$ be denoted by $x + m_fy = c_f$. If the bisectrix intersects $N(S)$ in the interior of a compact
edge, we denote the equation of this edge by $x + m_Sy = c_S$. Otherwise, let $x + m_Sy = c_S$ denote any line $l$ with 
$m_S$ and $c_S$ rational such that $N(S) \cap l = (d(S),d(S))$. Switching the roles of the $x$ and $y$ axes if necessary, we may
assume $m_S \leq m_f$. 

Let $r$ be a zero of $f_{e_0}(1,y)$ of order at least $d(N)$. We do a coordinate change now, letting $S_1(x,y) = S(x,y + rx^{m_f})$
and $f_1(x,y) = f(x,y +  rx^{m_f})$. (This is how the fractional powers of $x$ come in; we need that $m_S \leq m_f$ and
$m_f$ need not be an integer).

First we show that $d(S_1) \leq d(S)$. To see this, let $w$ denote the edge or vertex of $N(S)$ such that the supporting line of $N(S)$ of slope $-{1 \over m_f}$ intersects $N(S)$ at $w$, and let $z$ denote the edge or vertex of $N(S_1)$ such that the supporting line of $N(S_1)$ of slope $-{1 \over m_f}$ intersects $N(S_1)$ at $z$. Since $m_S \leq m_f$, either $w$ is a vertex
of $N(S)$ or $w$ is an edge of $N(S)$ intersecting or below the bisectrix. If $w$ is an edge or vertex lying entirely below the bisectrix then the coordinate change doesn't affect any of $N(S)$ on the bisectrix or higher, and thus $d(S_1) = d(S)$. If $w$ is a vertex of
$N(S)$ on or above the bisectrix, then the coordinate change results in an edge $z$ extending from $w$ to the $x$-axis. So $d(S_1) \leq d(S)$ here.  If $w$ is an edge
whose upper vertex is $(d(S),d(S))$, the coordinate change converts $w$ to an edge $z$ whose upper vertex is also $(d(S),d(S))$.
($z$ is an edge and not a vertex here since $S$ is in superadapted coordinates.) Thus $d(S_1) = d(S)$. Lastly, if 
$w$ is an edge of $N(S)$ intersecting the bisectrix in its interior then
 $(S_1)_{z}(1,y) = S_w(1, y+ r)$. Since $S$ is in superadapted coordinates, $S_w(1,y)$ has no zeroes of order $d(S)$ or higher (including at $y = 0$), so the same is true for $S_w(1,y + r)$. Thus the bisectrix intersects $N(S_1)$ at $(d(S),d(S))$, which is
in the interior of $z$. Hence $d(S_1) = d(S)$ here as well.

Note that the above argument shows that $(S_1)_z(x,y)$ necessarily has at least one term whose exponents lie below the bisectrix.
As a consequence, all supporting lines of $N(S_1)$ of slope
greater than $-{1 \over m_f}$ intersect $N(S_1)$ below the bisectrix only, a fact that will come in handy later on. 

We now examine the effect of the coordinate change on $N(f)$. Since $f$ has an edge $e_0$ of slope $-{1 \over m_f}$, $f_1$ will have an edge or vertex which we call $e'$ that is on the same line as $e_0$, such that $(f_1)_{e'}(1, y) = f_{e_0}(1, y + r)$. So since $f_{e_0}$ has a zero of
order $d(f)$ or greater at $y = r$, $(f_1)_{e'}(1, y)$ has a zero of order $d(f)$ or greater at $y = 0$. Hence $e'$ lies entirely
on or above the bisectrix. For any smooth function $g(x,y)$, it is 
not hard to show (as is shown in [G2] and elsewhere) that the maximum possible number of roots of any $g_e(x,y)$, counted 
according to multiplicity, is $2d(g)$. Hence in the case at hand, the maximum total number of roots of
$f_{e_0}(1,y)$ is $2d(f)$. Since $f_{e_0}(1,y)$ has a zero of order $d(f)$ or greater at $y = r$, the maximum possible order of  a zero of $(f_1)_{e'}(1, y) = f_{e_0}(1, y + r)$ besides $y = 0$ is $2d(f) - d(f) = d(f)$, with equality possible only if the bottom endpoint of $e'$ is on the bisectrix. 

 Let $(a,b)$ denote the lower vertex of $e'$, and let $\infty \geq m > m_f$ be minimal such that the line $l$
containing $(a,b)$ of slope $-{1 \over m}$ is horizontal or intersects either $N(f_1) - \{(a,b)\}$ or $N(S_1)$. If $m = \infty$, then
$N(S_1) \subset N(f_1)$ and the bisectrix intersects $N(f_1)$ in the horizontal ray, possibly at the 
vertex of this ray. In this case the conditions of Lemma 3.3 hold for the following reasons. If $e$ is an edge of $N(S_1)$ above $l$ with slope less than $-{1 \over m_f}$, the polynomials $S_1^{e}(1,y)$ and $f_1^{e}(1,y)$ are the same as
$S^{e}(1,y)$ and $f^{e}(1,y)$ since the coordinate change did not affect the terms of the polynomials on edges
above $e_0$.  So the terms of $S^{e}(1,y)$ and $f^{e}(1,y)$ have integral coordinates and lie wholly on or above the bisectrix. Thus
their zeroes are of order less than $d(f_1)$ with equality possible only if the lower endpoint of $e$ is on the bisectrix. If $N(f_1)$ has an edge of slope $-{1 \over m_f}$, it will be the edge we called $e'$ before. Recall that $(f_1)_{e'}(1,y) = f_{e_0}(1,y+ r)$, whose zeroes
we saw were of order at most $d(f) \leq d(f_1)$ where equality can hold here only if $(a,b)$ is a vertex of $N(f_1)$. We conclude
that we may apply Lemma 3.3. $f$ is a not a good direction here; by
 above $N(S_1)$ has a vertex $(v_1,v_2)$ below the bisectrix, so $v_2 < d(S_1)$. Since $N(S_1) \subset N(f_1)$ this vertex lies on or above the horizontal ray of $N(f_1)$. Hence for the generic $t$, $d(S_1 + tf_1) \leq d(f_1) \leq v_2 < d(S_1)$ and $f$ is not a good direction.

So from now on assume $m$ is finite. Except for finitely many
$t$, $N(S_1 + tf_1)$ will have an edge $k$ on the line $l$, and $(S_1 + tf_1)_k(x,y) = (S_1)^k(x,y) + t(f_1)^k(x,y)$. Let
$N_1$ denote the Newton polygon of $S_1 + tf_1$ for all but finitely many $t$, and $d(N_1)$ its Newton distance. 
We first consider the case where $(S_1)^k(x,y)$ is nonzero. By above all supporting lines of $N(S_1)$ of slope
greater than $-{1 \over m_f}$ intersect $N(S_1)$ entirely below the bisectrix. So since $-{1 \over m} > -{1 \over m_f}$, $(S_1)^k(x,y)$ is 
either $(S_1)_e(x,y)$ for some edge $e$ lying wholly below the bisectrix, or is of the form $cx^vy^w$ for some vertex 
$(v,w)$ of $N(S)$ lying below the bisectrix. Thus in either case, since $(a,b)$ is on or above the bisectrix, the bisectrix intersects $N_1$ in the interior of $k$ or at the upper vertex of $k$. 

 We are now in a position to apply Lemma 3.3. The polynomials
$(S_1)^{e}(1,y)$ and $(f_1)^{e}(1,y)$ for $e$  below $k$ are all below the bisectrix, so are all of degree less than $d(N_1)$ and thus have zeroes of
order less than $d(N_1)$.  The same is true for $(S_1)^k(1,y)$
since the exponents of the terms of $(S_1)^k(x,y)$ are all below the bisectrix. For any edge $e$ of $N_1$ above $k$ with slope $-{1 \over m_f}$ or smaller, the polynomials $S_1^{e}(1,y)$ and $f_1^{e}(1,y)$ have zeroes of small enough order to 
apply Lemma 3.3 similarly to two paragraphs ago.
We conclude that we may apply Lemma 3.3 to $S_1$ and $f_1$. Note that because $k$ 
intersects $N(S_1)$ below the bisectrix, one has $d(S_1 + tf_1) < d(S_1)$ for the generic $t$ and thus $f$ is not a good direction.

Thus it remains to consider the situation where $S_1^k(x,y)$ is zero. In this case $l$ intersects $N(f_1) - \{(a,b)\}$ but not 
$N(S_1)$. Thus the edge $k$ is also an edge of $N(f_1)$. If the lower vertex $(a',b')$ of $k$ is on or above the bisectrix,
we replace $(a,b)$ by $(a',b')$ and then repeat the argument of the previous three paragraphs. Since the $y$-coordinates
are integers, this can only happen finitely many times and eventually $(a',b')$ will be below the bisectrix. Denote by $k'$ the
analogue of $k$ once $(a',b')$ is below the bisectrix. Then if $k'$ intersects the bisectrix at its upper vertex, we may immediately apply Lemma 3.3 and we are done. Since this upper vertex is not on $N(S)$, $f$ is not a good direction in this situation. 
If $k'$ does not intersect the bisectrix at its upper vertex,
we are back to the situation of the sixth  paragraph of this proof, where we had $S^{e_0}(x,y)$ being the zero polynomial,
with the bisectrix intersecting $N$ in the interior of a compact edge $e_0$ coming entirely from $f$, so that $(S + tf)_{e_0}(x,y) = tf_{e_0}(x,y)$. The only difference is that instead of dealing with smooth functions we are now dealing with smooth functions of $x^{{1 \over \alpha}}$ and $y$, where $\alpha$ is the  difference between the $y$ coordinates of the upper and lower vertices of $e_0$.

Thus we may iterate the above process. If the procedure ends after finitely many iterations, we eventually
are in one of the cases already handled and we are done. Suppose the now that procedure does not ever end. Let $f_l$, 
$S_l$, and the Newton polygon $N_l$ correspond to the $l$th stage of the iteration, and let $d(N_l)$ denote the Newton distance. Let $(x,y)
\rightarrow (x, y + r_l x^{m_l})$ denote the coordinate change going from the $l$th to $l+1$st stage, and let
$(a_l,b_l)$ be the analogue of the lower vertex $(a,b)$ above. Since $b_l$ is nonincreasing with $l$ and has integral
values, eventually $b_l$ stops changing at some $l = l_0$ and stays at some value $b_0$. This means $a_l$ stays at some value $a_0$ as well. Furthermore, since the $b_l$ are nonincreasing, we 
are dealing with smooth functions of $x^{{1 \over b!}}$ and $y$ throughout. As a result, the slopes of the edges intersecting the bisectrix, increasing with each
iteration, will all be of the form ${-{nb! \over p }}$ with $n$ and $p$ integers such that $1 \leq n \leq b$. Because 
these slopes increase with each iteration and $n$ is constrained to be between $1$ and $b$, these slopes go to zero as 
$l$ goes to infinity. Since the vertex $(v_1,v_2)$ of $N(S_1)$ below the bisectrix will still be there in every iteration,
$b_0$ satisfies $b_0 \leq v_2 < d(S_1)$. 

We now use Borel's theorem (see [H]) to find an $s(x)$ which is a smooth function of $x^{{1 \over b!}}$ such that the Taylor expansion of $s(x)$ at $x = 0$ is $\sum_{l = l_0}^{\infty}r_lx^{m_l}$.
We perform the coordinate change $(x,y) \rightarrow (x,y + s(x))$. Define $\tilde{f}(x,y) = f_{l_0}(x,y + s(x))$ and 
$\tilde{S}(x,y) = S_{l_0}(x,y + s(x))$. Then $N(\tilde{f})$ has an infinite ray with vertex $(a_0,b_0)$,  $N(\tilde{S}) \subset N(\tilde{f})$, and $d(\tilde{S}) = d(S_1)$.  Now we may apply Lemma 3.3. Since $d(\tilde{S} + t\tilde{f}) \leq
d(\tilde{f}) = b_0 < d(S_1) = d(\tilde{S})$, $f$ is not a good direction here. 

We have now exhausted all the ways in which the algorithm may proceed. It is worth mentioning that in order to apply Lemma 3.3 as the iteration proceeds, we need first that the leftmost vertices of the Newton polygons arising have integral coordinates, which is automatic given the definition of the coordinate changes here, and secondly that the polynomials
coming from the new edges still satisfy the condition on the order of their zeroes. To see this, first observe that the edges not 
entirely on or above the bisectrix can be dealt with like before; the degree of any $S_l^e(x,y)$ 
and $f_l^e(x,y)$ appearing
will be at most $d(N_l)$, with equality only possible if the upper vertex of $e$ is on the bisectrix. Thus
the conditions needed for Lemma 3.3 will hold. For an edge $e'$ on or above the bisectrix, write its upper and lower endpoints 
as $(a'',b'')$ and $(a''',b''')$ respectively, and Then $b'' \leq b$, where $b$ is as before, and we saw $b \leq 2d(f) < 2d(N_l)$. Since $(a''',b''')$ is on or above the bisectrix,
we also have $b''' \geq d(N_l)$. Hence $b'' - b''' < d(N_l)$. This means any zero $f_{e'}(1,y)$ other than $y = 0$ is of order 
at most $b'' - b''' < d(N_l)$. Thus we may apply Lemma 3.3.

Thus we are done other than listing the good directions; that is, the $f$ for which $S + tf$ is of the same type as $S$ for
all but finitely many $t$. The only $f$ that were good directions in the above argument came near the beginning. If the bisectrix 
intersected $N(S)$ at a vertex or inside one of the rays, this only happened in the first paragraph of the main
argument in the case when $d(f) = d(N)$.
These correspond to the Newton polygon configurations given in Theorem 2.1. If the bisectrix intersected $N(S)$ inside the 
interior of a compact edge $e$, this only happened in the third and fourth paragraphs of the main argument, when $N(f)$ lay
 wholly on or above the line containing $e$, which is the Newton polygon configuration stipulated by Theorem 2.1 as well. This completes the proof. 

\noindent {\bf 4. Proof of Theorem 2.2.}

We first do the case where $S(x,y)$ is in Case 1 superadapted coordinates. Let $e$ denote
the edge of $N(S)$ intersecting the bisectrix in its interior, and let $x + my = c$ denote the equation of this edge. We focus
our attention on $(x,y)$ inside a small rectangle $T_r = [-r,r] \times [r^m,r^m]$, where $r$ is viewed as a fixed number
sufficiently small for our arguments to proceed. We do a finite Taylor expansion of $S(x,y)$ in the $x$ variable, obtaining
$$S(x,y) = \sum_{0 \leq j \leq c} {\partial_x^jS(0,y) \over j!}x^j + E(x,y) \eqno (4.1)$$
We then do a partial Taylor expansion in the $y$ variable of a given term ${\partial_x^jS(0,y) \over j!}$ up to the $y^k$ power, where
$k$ is maximal such that $j + mk \leq c$. Then $(4.1)$ becomes
$$S(x,y) = \sum_{0 \leq j + mk \leq c} {\partial_x^j\partial_y^k S(0,0) \over j!k!}x^jy^k + F(x,y) \eqno (4.2)$$
Here the error term $F(x,y)$ satisfies the estimate
$$|F(x,y)| \leq C(|x|^{c+1} + |y|^{{c \over m} + 1})\|S\|_{C^n(T_r)} \eqno (4.3)$$
Here $C$ is a constant depending only on the edge $e$ and $n = \max({c \over m}, c) + 1$, although the exact value of $n$ is
not important for us. Since the edge $e$ of $N(S)$ has equation $x + my = c$, $(4.2)$ can be rewritten as
$$S(x,y) = S_e(x,y) + F(x,y) \eqno (4.4)$$
Let $d = d(S)$. Suppose $1 \leq l < d$. Then $\partial_lS(x,y)$ has an edge or vertex $e_l$ on the line $x + my = c - ml$, and as in $(4.4)$ 
we may write
$$\partial_l^yS(x,y) = (\partial_l^y S)_{e_l}(x,y) + F_l(x,y) \eqno (4.5)$$
This time the error estimate is
$$|F_l(x,y)| \leq C(|x|^{c+1 - {lm}} + |y|^{{c \over m} - l + 1})\|S\|_{C^{n+l}(T_r)} \eqno (4.6)$$

Next, we divide $T_r$ into four parts along the curves $y = \pm|x|^m$. The estimates for each of the parts
is done the same way, so we focus on the $0 < x < r$, $-x^m < y < x^m$ region, which we denote by $U_r$. 
Note that on $U_r$, since $|y| < x^m$, the error estimates $(4.4)$ and $(4.6)$ become
$$|F(x,y)| \leq Cx^{c+\min(1,m)}\|S\|_{C^n(T_r)} \eqno (4.7a)$$
$$|F_l(x,y)| \leq Cx^{c - lm + \min(1,m)}\|S\|_{C^{n+l}(T_r)} \eqno (4.7b)$$
Suppose $1 \leq l < d$ is such that $S_e(1,y)$ has a zero of order $l$ at some $y = y_0$. Then since $S_e(x,y)$ is a mixed
homogeneous function, $S_e(x,y)$
has a zero of order $l$ in the $y$ direction at every point on the curve $y = y_0x^m$ for $x > 0$, and there is some $\epsilon
> 0$ such that on the wedge $W_{y_0,\epsilon,r} =\{ (x,y) \in U_r: (y_0 - \epsilon)x^m \leq y \leq (y_0 + \epsilon)x^m\}$ one has estimates
$$C_1x^{c - lm} < (\partial_l^y S)_{e_l}(x,y) < C_2x^{c - lm} \eqno (4.8)$$
Given the error estimate $(4.7b)$, if $x$ is sufficiently small, depending on $\|S\|_{C^{n+l}(T_r)}$, on $W_{y_0,\epsilon,r}$
 one has estimates
$$C_1'x^{c - lm} < (\partial_l^y S)(x,y) < C_2'x^{c - lm} \eqno (4.9)$$
Suppose now $0 < |y_0| < 1$ is such that $S_e(x,y)$ does not have a zero at $y_0$. Then along the curve $y = y_0x^m$, by mixed homogeneity $S_e(x,y)$ is 
equal to $ax^c$ for some nonzero $a$. Thus $\partial_x [S_e(x,y_0x^m)] = acx^{c-1}$ is nonzero. But by the chain rule
$$\partial_x [S_e(x,y_0x^m)] =  (\partial_x S)_{e^1}(x,y_0x^m) + my_0x^{m-1}(\partial_y S)_{e_1}(x,y_0x^m) \eqno (4.10)$$
Here $(\partial_x S)_{e^1}$ defined similarly to $(\partial_y S)_{e_1}$ except with the roles of the $x$ and $y$ variables reversed. Thus always
either $(\partial_x S)_{e^1}(x,y_0x^m)$ or $(\partial_y S)_{e_1}(x,y_0x^m)$ or nonzero. Therefore, by mixed homogeneity of $(\partial_x S)_{e^1}$ and  $(\partial_y S)_{e_1}$, either  $(\partial_x S)_{e^1}(x_0,1)  \neq 0$ for $x_0 = {y_0}^{1 \over m}$ or $(\partial_y S)_{e_1}(1,y_0) \neq 0$. In the latter case we have some $W_{y_0,\epsilon,r}$ such that 
$(4.9)$ holds with $l = 1$, and in the former case, reversing the roles of the $x$ and $y$ axes there is some 
analogous $W_{x_0,\epsilon,r}$ on which we have
$$C_1x^{c - 1} < \partial_x S(x,y) < C_2x^{c - 1} \eqno (4.11)$$
In the case where $y_0 = 0$, then we do a coordinate change $(x,y) \rightarrow (x, y + {1 \over 2}x^m)$. Then $y_0$ becomes
${1 \over 2}$ instead of 0, and now either $(4.9)$ holds for $l = 1$ or $(4.11)$ holds. The error estimates still hold and the possible introduction of fractional powers of $x$ does not interfere with our arguments.

Thus for any $y$, there is a small interval $[y_0 - \epsilon, y_0 + \epsilon]$ or $[x_0 - \epsilon, x_0 + \epsilon]$ with $x_0 = y_0^{1 \over m}$ such that on the associated wedge $W_{y_0,\epsilon,r}$  or $W_{x_0,\epsilon,r}$ either $(4.9)$ holds for some $1 \leq l < d$, or 
$(4.11)$ holds. (We use the fact that $S(x,y)$ is in 
superadapted coordinates here to ensure $l$ never has to be $d$ or higher). By compactness, we can cover
$T_r$ with finitely many such wedges $W_i$. 

Suppose now that $f(x,y)$ is a perturbation function whose Newton polygon lies entirely on or above the line $x + my = c$.
Then $(4.4)-(4.7)$ hold for $f$ in place of $S$. Hence there are constants $C_3$ and $C_4$ independent of $f$ such that if $\|f\|_{C^{n+d}}$ is sufficiently small, on a given $W_i$ one of the following two equations holds.
$$C_3x^{c - lm} < \partial_l^y(S + f) < C_4x^{c - lm} \eqno (4.12a)$$
$$C_3x^{c - 1} < \partial_x (S + f) < C_4x^{c - 1} \eqno (4.12b)$$
We now proceed as follows. If $(4.12a)$ holds, we perform a Van der Corput-type argument in the $y$ direction and integrate the result with respect to $x$, and if $(4.12b)$ holds we perform it in the $x$ direction
and integrate the result with respect to $y$. In the former case, Lemma 2.0 of [G2] gives
$$ |\int_{W_i}e^{i\lambda S(x,y)}\phi(x,y)\,dx\,dy| \leq C \lambda^{-{1 + m \over (c - lm) + ml}} \eqno (4.13a)$$
Since $d = {c \over 1 + m}$ , $(4.13a)$ is the same as 
$$ |\int_{W_i}e^{i\lambda S(x,y)}\phi(x,y)\,dx\,dy| \leq C \lambda^{-{1 \over d}} \eqno (4.13b)$$
In the proof of Lemma 2.0, the constant $C$ is of the form $\|\phi\|_{C^1(W_i)}$ times a function of $C_3$ and the supremum of finitely many derivatives of the phase $S + f$. In particular, we can write it as $C'\|\phi\|_{C^1(W_i)}$ where $C'$ is uniform
over all perturbations in question.  

\noindent In the case where $(4.12b)$ holds, one applies Lemma 2.0 in the $x$ direction,
obtaining 
$$ |\int_{W_i}e^{i\lambda S(x,y)}\phi(x,y)\,dx\,dy| \leq C \lambda^{-{1 + {1 \over m} \over ({c \over m} - {1 \over m}) + {1 \over m}}} \eqno (4.14)$$
The exponent on the right-hand side here is once again $-{m + 1 \over c} = -{1 \over d}$, so $(4.13b)$ holds once again with
$C$ of the form $C'\|\phi\|_{C^1(W_i)}$ This completes the proof of Theorem 2.2 when $S(x,y)$ is in Case 1 
superadapted coordinates. 

We now move to the case when $S(x,y)$ is in Case 2 superadapted coordinates. In this situation, there is a finite collection of
pairs $(l_1,l_2)$ of lines containing $(d,d)$ such that the Newton polygon of any good direction $f(x,y)$ is contained in
the wedge above $(d,d)$ determined by $l_1$ and $l_2$ for one of these pairs. It suffices to prove the uniform estimates
for a given pair. Write the equations of the lines $l_1$ and $l_2$ as $x + m_1y = c_1$ and $x + m_2y = c_2$. Let $m_3$  be any
rational number between $m_1$ and $m_2$.  Let $l_3$ be the line of slope $-{1 \over m_3}$ containing
$(d,d)$, and write its equation as $x + m_3y = c_3$. Thus $N(S)$ intersects $l_3$ at the one point $(d,d)$. This time,
we let $V_r$ be the rectangle $[-r,r] \times [-r^{m_3},r^{m_3}]$, and divide it into four parts via the curves $y = \pm |x|^{m_3}$.
We focus on the part where $0 < x < r$ and $-x^{m_3} < y < x^{m_3}$, which we denote by $W_r$, as the other three regions
are dealt with the same way.

We examine the expansion $(4.2)$ with $m$ and $c$ replaced by $m_3$ and $c_3$ respectively. For some
 $a \neq 0$ we get
$$S(x,y) = a x^dy^d + F(x,y) \eqno (4.15)$$
Here $|F(x,y)| \leq C|x|^{c_3 + \min(1,m_3)}\|S\|_{C^n(V_r)}$ as in $(4.7a)$. Similarly, taking the $d$th $y$ derivative gives 
$$\partial_y^d S(x,y) = a d! x^d + F_l(x,y) \eqno (4.16)$$
This time we have
$$|F_l(x,y)| \leq Cx^{c_3 - dm_3 + \min(1,m_3)}\|S\|_{C^{n+d}(V_r)} \eqno (4.17a)$$
Since $(d,d)$ is on $l_3$, one has $c_3 - dm_3 = d$, and the above becomes
$$|F_l(x,y)| \leq Cx^{d+ \min(1,m_3)}\|S\|_{C^{n+d}(V_r)} \eqno (4.17b)$$
Thus if $r$ is sufficiently small, depending on $\|S\|_{C^{n+d}(V_r)}$, we have
$$|\partial_y^d S(x,y)| > {|a| \over 2}x^d \eqno (4.18)$$
If the perturbation function $f$ is such that $|f|_{r,n + d}$ is sufficiently small, then looking at $(4.16)-(4.17)$ with
$f$ in place of $S$ gives
$$|\partial_y^d f(x,y)| < {|a| \over 4}x^d \eqno (4.19)$$
Combining with $(4.18)$ leads to
$$|\partial_y^d (S + f)(x,y)| > {|a| \over 4}x^d \eqno (4.20)$$
If we apply the Van der Corput lemma in the $y$ direction, and integrate the result in $x$ (see the proof of Lemma 4.1 of [G2] for
a very similar calculation), one obtains
$$ |\int_{W_r}e^{i\lambda S(x,y)}\phi(x,y)\,dx\,dy| \leq C \|\phi\|_{C_1(V_r)} \lambda^{-{1 + m \over dm + d}}\ln\lambda \eqno (4.21)$$
Since ${-{1 + m \over dm + d}}$ is just ${-{1 \over d}}$, $(4.21)$ is the desired estimate and we are done with the proof of
Theorem 2.2 when $S(x,y)$ is in Case 2 superadapted coordinates. 

We now move to the case where $S(x,y)$ is in Case 3 superadapted coordinates. Without loss of generality we may assume the
bisectrix intersects $N(S)$ in the interior of a horizontal ray. We denote the lower left vertex of $N(S)$ by $(c,d)$, where 
$0 \leq c < d = d(S)$, and we let $ax^cy^d$ be the corresponding term of the Taylor series of $S(x,y)$. If $r$ is sufficiently 
small, which we may assume, then on $D_r$ we have
$$|\partial_x^cS(x,y)| > c!{|a| \over 2}|y|^d \eqno (4.22a)$$
For $c = 0$ we may also assume that $r$ is small enough so that we have
$$|\partial_y^dS(x,y)| > d!{|a| \over 2} \eqno (4.22b)$$
The only conditions we need on $r$ is that $(4.22a)-(4.22b)$ holds, so in what follows we always assume $r$ is some fixed number
such that $(4.22a)-(4.22b)$ are satisfied.
Suppose $f(x,y)$ is an admissible direction. Then $f(x,y)$ is of the form $F(x,y)y^d$ for some real-analytic $F(x,y)$. The norm
$\|.\|_{r,N}$ is such that if $\eta$ is sufficiently small and $N$ is sufficiently large, then if $\|f\|_{r,N} < \eta$ we have
$$|\partial_x^cF(x,y)| < c!{|a| \over 4} \eqno (4.23a)$$
If $c = 0$ then if $\eta$ is sufficently small we also must have
$$|\partial_y^d (F(x,y)y^d)| < d!{|a| \over 4} \eqno (4.23b)$$
Thus in view of $(4.22)$, we have
$$|\partial_x^c(S + f)(x,y)| > c!{|a| \over 4}|y|^d \eqno (4.24a)$$
While if $c = 0$ we also have 
$$|\partial_y^d (S + f)(x,y)| > d! {|a| \over 4} \eqno (4.24b)$$
We use the Van der Corput lemma in conjunction with $(4.24a)-(4.24b)$ to get the desired upper bounds on $|J_{S + f,\phi}(\lambda)|$. When $c = 0$ we use $(4.24b)$, obtaining
$$|\int e^{i\lambda (S(x,y) + f(x,y))}\phi(x,y)\,dy| \leq C\lambda^{-{1 \over d}} \eqno (4.25)$$
Integrating this with respect to $x$ gives the desired uniform bounds. Suppose now $c > 0$. Then the Van der Corput lemma
applied to $(4.24a)$ gives 
$$|\int e^{i\lambda (S(x,y) + f(x,y))}\phi(x,y)\,dx| \leq C\max(1, \lambda^{-{1 \over c}} |y|^{-{d \over c}}) \eqno (4.26)$$
As a result, 
$$|J_{S + f,\phi}(\lambda)| \leq C \int_0^1 \max(1, \lambda^{-{1 \over c}}y^{-{d \over c}})\,dy \eqno (4.27)$$
The quantity $\lambda^{-{1 \over c}}y^{-{d \over c}}$ is equal to 1 when $y = \lambda^{-{1 \over d}}$. Since the 
exponent ${d \over c}$ is greater than 1, the integrals over the two sides of this breaking point will be comparable. The integral
over $y < \lambda^{-{1 \over d}}$ is just $\lambda^{-{1 \over d}}$, so we conclude that
$$|J_{S + f,\phi}(\lambda)|  \leq C' \lambda^{-{1 \over d}} \eqno (4.28)$$
These are the desired uniform estimates and we are done.

\noindent {\bf 5. Proof of Theorem 2.3.}

We first consider the case where $S(x,y)$ is in Case 1 superadapted coordinates. Since $|A_{\phi, S+f}|$ is uniformly bounded 
for $|f|_{r,N} < \eta$ by
Theorem 2.2, it suffices to show that there is some $\delta_ 0 > 0$ such that if $f_1$ and $f_2$ satisfies the hypotheses of
Theorem 2.3 and $|f_2 - f_1|_{r,N} < \delta_0$, then the estimate $(2.8)$ holds. This is what we will prove.

By $(A.3)$, $A_{\phi, S+f}$ is a linear combination of $B_{\phi,S + f}$ and $B_{\phi,-S - f}$, 
so it suffices to prove $(2.8)$ with $B_{\phi,S +f_i}$ in place of $A_{\phi,S + f_i}$ ($B_{\phi,-S - f_i}$ uses
 the same argument with $S$ replaced by $-S$, and $f_i$ replaced by $-f_i$). Furthermore, by $[G2]$ we have
$$B_{\phi, S} = (m + 1)^{-1}
\phi(0,0) \int_{-\infty}^{\infty} (S_e^+(1,y)^{-{ 1 \over d}} + S_e^+(-1,y)^{-{ 1 \over d}})\,dy \eqno (5.1)$$
Here $d = d(S)$, $e$ denotes the edge of $N(S)$ intersecting the bisectrix in its interior, and $S_e(x,y)^+$ denotes $\max(S_e(x,y),0)$. 
The $S_e^+(1,y)^{-{ 1 \over d}}$ and $ S_e^+(-1,y)^{-{ 1 \over d}}$ integrals are dealt with the same way, so we focus our
attention on the $S_e^+(1,y)^{-{ 1 \over d}}$ integral. Thus if we define
$$D_S = \int_{-\infty}^{\infty} S_e^+(1,y)^{-{ 1 \over d}}\,dy \eqno (5.2)$$
Then it suffices to show that for some constants $B$ and $\alpha$, whenever $|f_2|_{r,N}$,  $| f_1|_{r,N} < \eta$  we have
$$|D_{S+f_2} - D_{S + f_1}| \leq B |f_2 - f_1|_{r,N}^{\alpha} \eqno (5.3)$$
Write $S(1,y) = \sum_{i=m}^n c_i y^i$. Since $S(x,y)$ is in superadapted coordinates, $n > d$ and $m < d$. 
We break $D_S = D_S^1 + D_S^2$, where $D_S^1$ is the integral over $|y| < 1$ and $D_S^2$ is the integral over $|y| \geq 1$. 
We change variables from $y$ to ${1 \over y}$ in $D_S^2$, and obtain
$$D_S^2 = \int_{|y| \leq 1} y^{-2}((\sum_{i=m}^n c_i y^{-i})^+)^{-{1 \over d}}\,dy\eqno (5.4)$$
$$=  \int_{|y| \leq 1} ((\sum_{i=m}^n c_i y^{2d - i})^+)^{-{1 \over d}}\,dy\eqno (5.5)$$
We compare this expression with 
$$D_S^1 = \int_{|y| \leq 1} ((\sum_{i=m}^n c_i y^i)^+)^{-{1 \over d}}\eqno (5.6)$$
The highest power of $y$ appearing in $(5.5)$ is $2d - m > d $, and the lowest power is $2d -n < d$. Furthermore, perturbations
of $S$ in the $|.|_{r,N}$ norm will correspond to perturbations of the polynomial $\sum_{i=m}^n c_i y^{2d - i}$ in the 
$|.|_{r,N}$ norm; $(S + f)_e(x,y)$ depends solely  on finitely many terms of the Taylor expansion of $S + f$ at the origin. Thus
finding bounds of the form $(5.3)$ for $|D_{S+ f_2}^2 - D_{S+ f_1}^2|$ are equivalent to finding bounds of this form for $|D_{S+ f_2}^1 - D_{S+ f_1}^1|$ , so we restrict our attention to bounding the latter.

Since any zero of $S_e(1,y)$ is of order less than $d$, we can divide $[-1,1] = \cup_{i=1}^KI_i$, where $I_i$ is a closed interval
on which for some $\epsilon > 0$ and some $0 \leq l < d$ depending on $i$ we have an estimate 
$$\bigg|{d^l S_e(1,y) \over dy^l}\bigg| > 2\epsilon \eqno (5.7)$$
For a given perturbation function $f(x,y)$, we let $f_e(x,y)$ denote the sum of terms of $f$'s Taylor expansion lying on $e$, where $f_e(x,y)$ can now
be zero. Shrinking the $\eta$ and raising the $N$ of Theorem 2.2 if necessary, we may assume that for all $f$ with $|f|_{N,r}
< \eta$ on $I_i$ we have 
$$\bigg|{d^l f_e(1,y) \over dy^l}\bigg| < \epsilon$$
As a result we have
$$\bigg|{d^l (S + f)_e (1,y) \over dy^l}\bigg| > \epsilon \eqno (5.8)$$
Thus if $l \geq 1$, one may apply the measure Van der Corput lemma (see [C]) in the $y$ direction to obtain that there is a constant $C$ 
depending only on $S$ such that for all $\lambda$ the following holds.
$$|\{y \in I_i: |(S + f)_e(1,y)| < \lambda\}| < C\lambda^{1 \over l}\eqno (5.9)$$
Next, we write $D_S^1 = \sum_{i=1}^K E_S^i$, where $E_S^i$ denotes the portion of the integral over $I_i$. We will give a bound 
of the form $(5.3)$ for each $|E_{S+f_2}^i - E_{S + f_1}^i|$. Note that
$$E_{S+f_2}^i - E_{S + f_1}^i \leq \int_{I_i}\big|((S + f_2)_e(1,y))^+)^{-{1 \over d}} -  ((S + f_1)_e(1,y))^+)^{-{1 \over d}}\big| \eqno (5.10)$$
Denote $|f_2 - f_1|_{r,N}$ by $\delta$. We write $(5.10)$ as $J_1 + J_2$,
where $J_1$ is the integral over $y$ for which
$|(S + f_1)_e(1,y)| < \delta^a$ for some small $a$ to be determined by our arguments, and $J_2$ is the integral over $y$
where $|(S + f_1)_e(1,y)| \geq \delta^a$. Note that if $l = 0$ in $(5.8)-(5.9)$, then if $\delta < \epsilon^{1 \over a}$, which we may assume, by $(5.8)$ $|(S + f_1)_e(1,y)|$ is never less than $\delta^a$ and therefore $J_1 = 0$. 
Thus when analyzing $J_1$ we may always assume $l \geq 1$. By $(5.10)$ we have
$$|J_1| \leq  \int_{|(S + f_1)_e(1,y)| < \delta^a} |(S + f_2)_e(1,y)|^{-{1 \over d}} +  \int_{|(S + f_1)_e(1,y)|< \delta^a}|(S + f_1)_e(1,y)|^{-{1 \over d}} \eqno (5.11)$$
Since $|f_2 - f_1|_{r,N} = \delta$, $|(f_2 - f_1)_e(1,y)| < C\delta$ for all $y \in I_i$. Therefore $|(S + f_2)_e(1,y)|
\leq |(S + f_1)_e(1,y)| + |(f_1 - f_2)_e(1,y)|  < \delta^a + C\delta < 2\delta^a$. Hence for any $y$ in the 
domain of the first integral, we have $|(S + f_2)_e(1,y)| < \delta + \delta^a < 2\delta^a$. Thus $(5.11)$ is bounded by
$$|J_1| \leq  \int_{|(S + f_2)_e(1,y)| < 2\delta^a} |(S + f_2)_e(1,y)|^{-{1 \over d}} +  \int_{|(S + f_1)_e(1,y)|< \delta^a}|(S + f_1)_e(1,y)|^{-{1 \over d}} \eqno (5.12)$$
We bound the first term of $(5.12)$ as the second term is done the same way. For any $g(y)$, by the relation between $L^p$ 
norms and distribution functions, applied to ${1 \over |g(y)|}$, we have
$$\int_{|g(y)| < \lambda_0}|g(y)|^{-{1 \over d}}\,dy = {1 \over d}\int_0^{\lambda_0}\lambda^{-{1 \over d} - 1}|\{y: |g(y)| < \lambda\}|\,d\lambda \eqno (5.13)$$
Thus the first term of $(5.12)$ is bounded by
$${1 \over d}\int_0^{2\delta^a}\lambda^{-{1 \over d} - 1}|\{y: |(S + f_2)_e(1,y)| < \lambda\}|\,d\lambda \eqno (5.14)$$
Inserting $(5.9)$ (recalling $l \geq 1$), we have that $(5.14)$ is bounded by
$${1 \over d}\int_0^{2\delta^a}\lambda^{-{1 \over d} + {1 \over l} - 1}\,d\lambda \eqno (5.15)$$
Since $l < d$, the exponent in $(5.15)$ is greater than $-1$, so $(5.15)$ is bounded by $C\delta^{a'}$ for some $a' > 0$.
This is the estimate that we seek since $|f_2 - f_1|_{r,N} = \delta$.

It remains to bound $J_2$, the portion of $(5.10)$ where $|(S + f_1)_e(1,y)| \geq \delta^a$. In this case, $|(S + f_2)_e(1,y) 
- (S + f_1)_e(1,y)|  = |(f_2 - f_1)_e(1,y)|  \leq C\delta < {1 \over 2}\delta^a$. So if the integrand in $J_2$ is nonzero the terms are both positive. Thus $J_2$ becomes 
$$\int_{ (S + f_1)_e(1,y) \geq  \delta^a} \big|((S + f_2)_e(1,y))^{-{1 \over d}} -  ((S + f_1)_e(1,y))^{-{1 \over d}}\big| \eqno (5.16)$$
By the mean value theorem, for each $y$ there is a $t_y$ with $|t_y| \leq   |(f_1 - f_2)_e(1,y)|  < C\delta$ such that
the integrand of $(5.16)$ is equal to $ {1 \over d}|(f_2 - f_1)_e(1,y)  ((S + f_1)_e(1,y) + t_y)^{-{1 \over d} - 1}|$.  Since
$|t_y| < C\delta$ and $(S + f_1)_e(1,y) \geq \delta^a$,  $|(S + f_1)_e(1,y) + t_y|$ is at least ${1  \over 2}(S + f_1)_e(1,y)$. So
since $|(f_2 - f_1)_e(1,y)| < C\delta$, we conclude that the integrand in $(5.16)$ is bounded by $C\delta (S + f_1)_e(1,y)^{-{1 \over d} - 1}$. Thus $(5.16)$ is bounded by
$$C\delta \int_{|(S + f_1)_e(1,y)| \geq  \delta^a}|(S + f_1)_e(1,y)|^{-{1 \over d} - 1}\,dy \eqno (5.17)$$
If $l = 0$, then the integrand in $(5.17)$ is bounded below by $(5.8)$ and thus $(5.17)$ gives a bound of $C\delta$, better than
what we need. If $l > 0$, we again use the relation between $L^p$ norms and distribution functions, this time in the form
$$\int_{(S + f_1)_e(1,y) \geq  \delta^a}|(S + f_1)_e(1,y)|^{-{1 \over d}- 1}\,dy = \big({1 + {1 \over d}}\big) \int_{\delta^a}^{\infty} \lambda^{-{1 \over d} - 2}|\{y: |(S + f_1)_e(1,y)| < \lambda\}|\,d\lambda \eqno (5.18)$$
Substituting $(5.9)$ in this time gives that $(5.17)$ is bounded by
$$C\delta \int_{\delta^a}^{\infty} \lambda^{{1 \over l} - {1 \over d} - 2}\,d\lambda \eqno (5.19)$$
Since $l \geq 1$, the exponent here is less than $-1$, and we may integrate. Thus $(5.19)$ becomes $C\delta^{1 + a({1 \over l}
-{1 \over d} - 1)}$. As long as $a$ is sufficiently small this is bounded by say $C\delta^{1 \over 2}$ and we are done with the
proof for the case of Case 1 superadapted coordinates. 

We now proceed to Case 3 superadapted coordinates in the real-analytic case. Without loss of generality, we assume that
the bisectrix intersects $N(S)$ in the interior of a horizontal ray. We let $a(x)y^d$ denote the sum of the terms of $S(x,y)$'s
Taylor expansion on this ray, so that $a(x)$ is real-analytic and $d = d(S)$. In [G2] it is proven that $(A.3)$ holds in the
real-analytic situation for Case 3 superadapted coordinates, so like above it suffices to prove $(2.9)$ with $B_{\phi,S}$ in
place of $A_{\phi,S}$. Explicit formulas for $B_{\phi,S}$ are given in [G2]. If $d$ is even, we have
$$B_{\phi,S} = 2\int_{-\infty}^{\infty} (a(x)^+)^{-{1 \over d}}\phi(x,0)\,dx \eqno (5.20a)$$
While if $d$ is odd, we have
$$B_{\phi,S} = \int_{-\infty}^{\infty} |a(x)|^{-{1 \over d}}\phi(x,0)\,dx \eqno (5.20b)$$
For an allowable perturbation function $f(x,y)$, we may write $f(x,y) = F(x)y^d+ O(y^{d+1})$, with $F(x)$ real-analytic,
 so that when $d$ is even we have
$$B_{\phi,S + f_2} - B_{\phi,S + f_1} = 2\int_{-\infty}^{\infty}\big(((a(x) + F_2(x))^+)^{-{1 \over d}} - ((a(x) + F_1(x))^+)^{-{1 \over d}}\big)\phi(x,0)\,dx \eqno (5.21a)$$
While if $d$ is odd we have
$$B_{\phi,S + f_2} - B_{\phi,S + f_1} = \int_{-\infty}^{\infty}\big(|(a(x) + F_2(x))|^{-{1 \over d}} - |(a(x) + F_1(x))|^{-{1 \over d}}\big)\phi(x,0)\,dx \eqno (5.21b)$$
Thus if $d$ is even we have
$$|B_{\phi,S + f_2} - B_{\phi,S + f_1}| \leq 2\|\phi\|_{L^{\infty}}\int_{-r}^{r}\big|((a(x) + F_2(x))^+)^{-{1 \over d}} - ((a(x) + F_1(x))^+)^{-{1 \over d}}\big|\,dx\eqno (5.22a)$$
And if $d$ is odd we have
$$|B_{\phi,S + f_2} - B_{\phi,S + f_1}| \leq \|\phi\|_{L^{\infty}}\int_{-r}^{r}\big||(a(x) + F_2(x))|^{-{1 \over d}} - |(a(x) + F_1(x))|^{-{1 \over d}}\big|\,dx \eqno (5.22b)$$
Recall $r$ is the radius of the disk we are working in. Let $ax^c$ denote the term of lowest order of $a(x)$. Then we may assume
that $r$ is sufficiently small that $|{d^c (a(x) + F(x)) \over dx^c}| > {|a| \over 2}$ on the disk for all allowable perturbations.
Thus we may argue as in Case 1 with just one $I_i$, using this inequality in place of $(5.8)$ and using the fact that
$|F_2(x) - F_1(x)| \leq ||f_2 - f_1||_{r,N} $. The argument for $(5.22a)$ is 
exactly the same as in Case 1, while in $(5.22b)$ one makes minor modifications due to the fact one is no longer taking the 
positive parts of the functions in question. We omit the details for brevity and we are done.

\noindent {\bf Appendix. Formulas for the coefficient of the principal term of the asymptotics.}

In this appendix, we show that the formulas for the coefficients $A_{S,\phi}$ given in the real-analytic case in [G2] carry
over to the smooth situation when $S(x,y)$ is in Case 1 or Case 2 superadapted coordinates. When $S(x,y)$ has a Morse
critical point at the origin, explicit formulas are well-known for the general smooth case (see [S] p.344-347), so throughout
we will assume $S(x,y)$ does not have a Morse critical point at the origin. In [G2], one defines
$$I_{S, \phi}(\epsilon) =  \int_{\{(x,y): 0 < S(x,y) < \epsilon\}} \phi(x,y) \,dx\,dy \eqno (A.1)$$
It is shown in [G2] that for Case 1 or Case 2 in the smooth case, and in Case 3 for the real-analytic case only, that if $S(x,y)$ is in
 superadapted coordinates then as $\epsilon \rightarrow 0$ one has asymptotics
$$I_{S,\phi}(\epsilon) = B_{\phi,S}\,\,\epsilon^{\delta}|\ln(\epsilon)|^p + o(\epsilon^{\delta}|\ln(\epsilon)|^p) \eqno (A.2)$$
Here $(-\delta,p)$ is the oscillatory type of $S$ and explicit formulas for $B_{\phi,S}$ are shown in [G2]. In the real-analytic case,
one can use well-known methods (cf Ch 7 of [AGV]) to get explicit formulas for the $A_{\phi,S}$ in terms of these formulas for
$B_{\phi,S}$. Namely, in Case 1 and Case 3 superadapted coordinates we have
$$ A_{\phi,S} = {\Gamma({1 \over d}) \over d}(e^{i {\pi \over 2d}} B_{\phi,S} + e^{-i{ \pi \over 2d}}
B_{\phi,-S}) \eqno (A.3)$$
In Case 2, one has
$$A_{\phi,S} = -{\Gamma({1 \over d}) \over d}(e^{i {\pi \over 2d}} B_{\phi,S}  + e^{-i{ \pi \over 2d}}
B_{\phi,-S}) \eqno (A.4)$$
In this section, we will explain why formulas $(A.3)$ and $(A.4)$ still hold for $S(x,y)$ if it is in Case 1 or Case 2 superadapted
coordinates. In both cases our arguments resemble those used in [IM2] for the Case 1 situation. (Case 2 is analyzed differently in that paper.) 

We first suppose $S(x,y)$ is in Case 1 superadapted coordinates. As in section 4, we denote by $e$ the edge of $N(S)$ 
intersecting the bisectrix in its interior, and the equation of the line containing this edge by $x + my = c$. Again we 
take a small enough rectangle $[-r,r] \times [-r^m,r^m]$ and divide it into 4 regions via the curves $y = \pm|x|^m$. As before 
we focus our attention on the the region $U_r$ consisting of points where $0 < x < r$ and $-x^m < y < x^m$ as the other
three regions are dealt with similarly. We also again use the wedges $W_i$ used in section 4. Recall each wedge $W_i$ is of the
form $\{(x,y): 0 < x < r,\, a_ix^m < y < b_ix^m\}$, and on a given wedge we either have an estimate $\partial_y S(x,y) > 
Cx^{c - lm}$ for some $1 \leq l < d(S)$, or an estimate $\partial_xS(x,y) > Cx^{c - 1}$ (possibly after a coordinate 
change of the form $(x,y) \rightarrow (x,y + {1 \over 2}x^m)$ in the case where the wedge was centered along the $x$-axis.)

Let $\epsilon > 0$ be a small number, to be determined by our arguments. We write $W_i = X_i \cup Y_i$, where
$$X_i = \{(x,y): 0 < x < \lambda^{-\epsilon},\,a_ix^m < y < b_ix^m\} \eqno (A.5a)$$
$$Y_i = \{(x,y): \lambda^{-\epsilon} < x < r,\,a_ix^m < y < b_ix^m\} \eqno (A.5b)$$
We first estimate $\int_{Y_i}e^{i\lambda S(x,y)}\phi(x,y)\,dx\,dy$. To do this, if  $\partial_y S(x,y) > Cx^{c - lm}$ on $W_i$
then we use the Van der Corput lemma for measures (see [C]) in the $y$ direction and integrate the result with respect to $x$, and if $\partial_xS(x,y) > Cx^{c - 1}$ we use the Van der Corput lemma for measures in the $x$ direction and integrate the result 
with respect to $y$. In the former case, for a given $x$ we get the estimate
$$|\int_{a_i x^m}^{b_ix^m}e^{i\lambda S(x,y)}\phi(x,y)\,dy| \leq C\lambda^{-{1 \over l}} |x|^{-{c \over l} + m}\eqno (A.6)$$
Thus we have
$$|\int_{Y_i}e^{i\lambda S(x,y)}\phi(x,y)\,dx\,dy| \leq C \lambda^{-{1 \over l}}\int_{\lambda^{-{\epsilon}}}^r
|x|^{-{c \over l} + m}\,dx \eqno (A.7)$$
Since $l < d$, we may let $\delta > 0$ be the minimal possible value of ${1 \over l} - {1 \over d(S)}$. Thus if we choose $\epsilon$ small enough that the $x$ integral here is bounded by $C\lambda^{{\delta \over 2}}$, then since $-{1 \over l} + {\delta \over 2}
\leq -{1 \over d(S)} - {\delta \over 2}$ we have 
$$|\int_{Y_i}e^{i\lambda S(x,y)}\phi(x,y)\,dx\,dy| \leq C' \lambda^{-{1 \over d(S)} - {\delta \over 2}} \eqno (A.8)$$
Hence the $Y_i$ integral is of order of magnitude less than $\lambda^{-{1 \over d(S)}}$. We now turn to the $X_i$ integral. Here we 
write $S(x,y) = S_N(x,y) + (S(x,y) - S_N(x,y))$, where $S_N(x,y)$ is the sum of the terms of the Taylor expansion of 
$S(x,y)$ to a high enough order $N$ to be determined as a function of $\epsilon$. Then 
$$e^{i\lambda S(x,y)} = e^{i\lambda S_N(x,y)}
e^{i\lambda (S(x,y) - S_N(x,y))} \eqno (A.9)$$
$$ = e^{i\lambda S_N(x,y)} +e^{i\lambda S_N(x,y)} (e^{i\lambda (S(x,y) - S_N(x,y))} - 1) \eqno (A.10)$$ Since
on $X_i$ we have $0 < x < \lambda^{-\epsilon}$ and $|y| < x^m$, by choosing $N$ sufficiently large (depending on $\epsilon$) we can ensure that on
$X_i$ we have 
$$|S(x,y) - S_N(x,y)| \leq C_N\lambda^{-3} \eqno (A.11)$$
By the mean value theorem this in turn implies that 
$$|e^{i\lambda S_N(x,y)} (e^{i\lambda (S(x,y) - S_N(x,y))} - 1)| \leq C_N'\lambda^{-2} \eqno (A.12)$$
Thus we have 
$$\int_{X_i}e^{i\lambda S(x,y)}\phi(x,y)\,dx\,dy = \int_{X_i}e^{i\lambda S_N(x,y)}\phi(x,y)\,dx\,dy + O(\lambda^{-2})\eqno (A.13)$$
As a result,
$$\int_{W_i}e^{i\lambda S(x,y)}\phi(x,y)\,dx\,dy = \int_{X_i}e^{i\lambda S_N(x,y)}\phi(x,y)\,dx\,dy + O(\lambda^{-{1 \over d(S)} - {\delta \over 2}}) \eqno (A.14)$$
Furthermore, 
by applying the $Y_i$ integral argument to to $S_N$ in place of $S$ for $N$ sufficiently large, we have
$$\int_{X_i}e^{i\lambda S_N(x,y)}\phi(x,y)\,dx\,dy= \int_{W_i}e^{i\lambda S_N(x,y)}\phi(x,y)\,dx\,dy + O( \lambda^{-{1 \over d(S)} - {\delta \over 2}}) \eqno (A.15)$$
 We conclude that 
$$\int_{W_i}e^{i\lambda S(x,y)}\phi(x,y)\,dx\,dy = 
\int_{W_i}e^{i\lambda S_N(x,y)}\phi(x,y)\,dx\,dy + O( \lambda^{-{1 \over d(S)} - {\delta \over 2}}) \eqno (A.16)$$
Adding this over all $W_i$ and over all 4 pieces of 
$[-r,r] \times [-r^m,r^m]$, we get that if $N$ is sufficiently large, then
$$J_{S, \phi}(\lambda) = J_{S_N,\phi}(\lambda) + O( \lambda^{-{1 \over d(S)} - {\delta \over 2}}) \eqno (A.17)$$
But $J_{S_N,\phi}(\lambda)$ has
asymptotics with an $O(\lambda^{-{1 \over d(S)}})$ term, whose formula is given explicitly in terms of finitely many terms of the Taylor expansion of $S_N(x,y)$. Since these terms are the same for $S$ and $S_N$ if $N$ is sufficiently large, these formulas will
hold for $S(x,y)$ as well. This completes the proof for Case 1 superadapted coordinates.

We now move on to Case 2 superadapted coordinates. We argue similarly to the Case 1 situation, using the Case 2 wedges of
Section 4. Specifically, we let $m$ be such that the line $x + my = c$ intersects $N(S)$ at the point $(d(S),d(S))$ only, and
we divide the rectangle $[-r,r] \times [-r^m, r^m]$ into 4 regions $V_r$ via the curves $y = \pm |x|^m$. As in section 4, we 
focus our attention on $W_r$, the region where $0 < x < r$ and $-r^m < y < r^m$ as the other three regions are done the
same way. As in $(4.18)$, on $W_r$ we have an estimate
$$|\partial_y^dS(x,y)| > Cx^d \eqno (A.18)$$
Similar to in Case 1, for small $\epsilon > 0$ we subdivide into regions $X$ and $Y$ defined by
$$X= \{(x,y): 0 < x < \lambda^{-\epsilon}, -x^m < y < x^m\} \eqno (A.19a)$$
$$Y= \{(x,y): \lambda^{-\epsilon} < x < r, -x^m < y < x^m\} \eqno (A.19b)$$
One may use the Van der Corput lemma in the $y$ direction to obtain
$$|\int_{-x^m}^{x^m}e^{i\lambda S(x,y)}\phi(x,y)\,dy| \leq C\lambda^{-{1 \over d}} |x|^{-1}\eqno (A.20)$$
Thus we have
$$|\int_{Y}e^{i\lambda S(x,y)}\phi(x,y)\,dx\,dy| \leq C \lambda^{-{1 \over d}}\int_{\lambda^{-{\epsilon}}}^r
|x|^{-1}\,dx \eqno (A.21)$$
$$\leq C\epsilon\ln(\lambda)\lambda^{-{1 \over d}}\eqno (A.22)$$
So while this does contribute to the main term of the asymptotics, it does so in a way that shrinks linearly with $\epsilon$.
As for the $X$ integral, one can argue exactly as in the Case 1 situation and say that if one replaces $S$ by $S_N$ in 
$\int_{X}e^{i\lambda S(x,y)}\phi(x,y)\,dx\,dy$, the difference is bounded in absolute value by $C_{\epsilon}\lambda^{-2}$. (Since $N$ is a 
function of $\epsilon$ so is the constant here.) Furthermore, the argument used to show $(A.22)$ works for $S_N$ in place of $S$, so $|\int_{Y} e^{i\lambda S_N(x,y)}\phi(x,y)\,dx\,dy|$ is also bounded by
$(A.22)$. We conclude that
$$\bigg|\int_{W_r}e^{i\lambda S(x,y)}\phi(x,y)\,dx\,dy - \int_{W_r}e^{i\lambda S_N(x,y)}\phi(x,y)\,dx\,dy\bigg| \leq C_{\epsilon}\lambda^{-2} +  C\epsilon\ln(\lambda)\lambda^{-{1 \over d}} \eqno (A.23)$$
Adding $(A.23)$ over all four squares, we conclude that
$$|J_{S,\phi}(\lambda) -  J_{S_N,\phi}(\lambda)| \leq C_{\epsilon}\lambda^{-2} +  C\epsilon\ln(\lambda)\lambda^{-{1 \over d}}\eqno (A.24)$$
Like in Case 1, the  formulas of [G2] are such that for $N$ sufficiently large, the formula applied to $S$ is the 
same as the formula applied to $S_N$. Thus we have shown that the formula for the leading term of the asymptotics in 
the real-analytic case also holds in the smooth case modulo a term bounded by 
$C\epsilon\ln(\lambda)\lambda^{-{1 \over d}}$ as $\lambda \rightarrow \infty$. Letting $\epsilon$ go to zero shows that the leading terms are in fact the same, and we are done.

\noindent {\bf References.}

\noindent [AGV] V. Arnold, S. Gusein-Zade, A. Varchenko, {\it Singularities of differentiable maps},
Volume II, Birkhauser, Basel, 1988. \parskip = 3pt\baselineskip = 3pt

\noindent [C] M. Christ, {\it Hilbert transforms along curves. I. Nilpotent groups}, Annals of 
Mathematics (2) {\bf 122} (1985), no.3, 575-596.

\noindent [DK] J-P Demailly, J. Koll\'ar, {\it Semi-continuity of complex singularity exponents and K\"ahler-Einstein metrics on Fano orbifolds}, Ann. Sci. \'Ecole Norm. Sup. {\bf 34} (2001), no. 4, 525-556.

\noindent [G1] M. Greenblatt, {\it Resolution of singularities in two dimensions and the stability of integrals}, Adv. Math., 
{\bf 226} no. 2 (2011) 1772-1802.

\noindent [G2] M. Greenblatt, {\it The asymptotic behavior of degenerate oscillatory integrals in two
dimensions}, J. Funct. Anal. {\bf 257} (2009), no. 6, 1759-1798.

\noindent [H] L. H$\ddot{{\rm o}}$rmander, {\it The analysis of linear partial differential operators. I. 
Distribution theory and Fourier analysis}, 2nd ed. Springer-Verlag, Berlin, (1990).
xii+440 pp. 

\noindent [K1] V. N. Karpushkin, {\it A theorem concerning uniform estimates of oscillatory integrals when
the phase is a function of two variables}, J. Soviet Math. {\bf 35} (1986), 2809-2826.

\noindent [K2] V. N. Karpushkin, {\it Uniform estimates of oscillatory integrals with parabolic or 
hyperbolic phases}, J. Soviet Math. {\bf 33} (1986), 1159-1188.

\noindent [IKeM] I. Ikromov, M. Kempe, and D. M\"uller, {\it Estimates for maximal functions associated
to hypersurfaces in $R^3$ and related problems of harmonic analysis}, Acta Math. {\bf 204} (2010), no. 2,
151--271.

\noindent [IM1] I. Ikromov, D. M\"uller, {\it On adapted coordinate systems}, Trans. AMS, {\bf 363} (2011), 2821-2848.

\noindent [IM2] I. Ikromov, D. M\"uller, {\it  Uniform estimates for the Fourier transform of surface-carried measures in
$\R^3$ and an application to Fourier restriction}, J. Fourier Anal. Appl, {\bf 17} (2011), no. 6, 1292-1332.

\noindent [PSSt] D. H. Phong, E. M. Stein, J. Sturm, {\it On the growth and stability of real-analytic 
functions}, Amer. J. Math. {\bf 121} (1999), no. 3, 519-554.

\noindent [S] E. Stein, {\it Harmonic analysis; real-variable methods, orthogonality, and oscillatory 
integrals}, Princeton Mathematics Series Vol. 43, Princeton University Press, Princeton, NJ, 1993.

\noindent [V] A. N. Varchenko, {\it Newton polyhedra and estimates of oscillatory integrals}, Functional 
Anal. Appl. {\bf 18} (1976), no. 3, 175-196.

\noindent 
\end